\documentclass[11pt, twoside]{amsart}

\usepackage[utf8]{inputenc}
\usepackage[T1]{fontenc}
\usepackage{amssymb,amsmath,amstext, graphicx}

\newtheorem{theor}{Theorem}[section]
\newtheorem{lem}[theor]{Lemma}
\newtheorem{defin}[theor]{Definition}

\newtheorem{notation}[theor]{Notation}
\newtheorem{exam}[theor]{Example}

\newtheorem{cor}[theor]{Corollary}
\newtheorem{rem}[theor]{Remark}

\newtheorem{fact}[theor]{Fact}

\newtheorem{assump}[theor]{Assumption}

\numberwithin{equation}{section}

\newcommand{\dom}{\mathrm{dom}}
\newcommand{\acl}{\mathrm{acl}}
\newcommand{\dcl}{\mathrm{dcl}}
\newcommand{\tp}{\mathrm{tp}}

\newcommand{\es}{\emptyset}

\newcommand{\su}{\mathrm{SU}}

\newcommand{\nts}{\negthickspace}
\newcommand{\uhrc}{\nts \upharpoonright \nts}

\newcommand{\meq}{^{\mathrm{eq}}}

\newcommand{\mcA}{\mathcal{A}}
\newcommand{\mcB}{\mathcal{B}}
\newcommand{\mcC}{\mathcal{C}}
\newcommand{\mcD}{\mathcal{D}}

\newcommand{\mcF}{\mathcal{F}}
\newcommand{\mcG}{\mathcal{G}}

\newcommand{\mcM}{\mathcal{M}}
\newcommand{\mcN}{\mathcal{N}}

\newcommand{\mcP}{\mathcal{P}}

\newcommand{\mbP}{\mathbf{P}}

\newcommand{\mbR}{\mathbf{R}}
\newcommand{\mbK}{\mathbf{K}}

\newcommand{\us}{\underset}
\newcommand{\ind}{\raisebox{-2pt}[5pt][0pt]{$\smile$} \hspace*{-6.8pt}\raisebox{3pt}[5pt][0pt]{$|$} \; \: }
\newcommand{\nind}{\raisebox{-2pt}[5pt][0pt]{$\smile$} 
\hspace*{-6.8pt}\raisebox{3pt}[5pt][0pt]{$|$}\hspace*{-6.8pt}
\raisebox{3pt}[5pt][0pt]{$\diagup$} }

\newcommand{\rng}{\mathrm{rng}}

\title[Simple homogeneous structures]
{On sets with rank one in simple homogeneous structures}

\author[O. Ahlman]{Ove Ahlman}
\address{Department of Mathematics, Uppsala University, Box 480,
75106 Uppsala, Sweden.}
\email{ove@math.uu.se}

\author[V. Koponen]{Vera Koponen}
\address{Department of Mathematics, Uppsala University, Box 480,
75106 Uppsala, Sweden.}
\email{vera@math.uu.se}

\begin{document}

\begin{abstract}
We study definable sets $D$ of SU-rank 1 in $\mcM\meq$, where $\mcM$ is a countable homogeneous and simple structure in
a language with finite relational vocabulary. Each such $D$ can be seen as a `canonically embedded structure', which inherits
all relations on $D$ which are definable in $\mcM\meq$, and has no other definable relations. 
Our results imply that if no relation symbol of the language of $\mcM$ has arity higher than 2, then
there is a close relationship between triviality of dependence and $\mcD$ being a reduct of a binary random structure.
Somewhat more precisely:
(a) if for every $n \geq 2$, every $n$-type $p(x_1, \ldots, x_n)$ which is realized in $D$ is determined by its sub-2-types $q(x_i, x_j) \subseteq p$,
then the algebraic closure restricted to $D$ is trivial;
(b) if $\mcM$ has trivial dependence, then $\mcD$ is a reduct of a binary random structure.
\end{abstract}



\maketitle

\section{Introduction}\label{Introduction}

\noindent
We call a countable first-order structure $\mcM$ {\em homogeneous} if it has a finite relational vocabulary (also called signature)
and every isomorphism between finite substructures of $\mcM$ can be extended to an automorphism of $\mcM$.
(The terminology  {\em ultrahomogeneous} is used in some texts.)
For surveys about homogeneous structures and connections to other areas, see \cite{Mac10} and
the first chapter of \cite{Che98}.
It is possible to construct $2^{\omega}$ countable homogeneous structures, even for a vocabulary with only a binary relation symbol,
as shown by Henson \cite{Hen72}. But it is also known that in several cases,
such as partial orders, undirected graphs, directed graphs or stable structures with finite relational vocabulary, 
all countable homogeneous structures in each class can be classified in a more or less explicit way
\cite{Che98, Gar, GK, JTS, Lach84, Lach97, LT, LW, Schm, Shee}.
Ideas from stability theory and the study of homogeneous structures have been used to obtain a good understanding
of structures that are $\omega$-categorical and $\omega$-stable (which  need not be homogeneous) \cite{CHL} and,
more generally, of smoothly approximable structures \cite{CH, KLM}.

Simplicity \cite{Cas, Wag} is a notion that is more general than stability.
The structures that are stable, countable and homogeneous are well understood, 
by the work of Lachlan and others; see for example the survey \cite{Lach97}.
However, little appears to be known about countable homogeneous structures that are simple, even for a {\em binary vocabulary},
i.e. a {\em finite} relational vocabulary where every relation symbol has arity at most 2. 
Besides the present work, \cite{Kop-one-based} and the dissertation of Aranda L\'{o}pez \cite{AL} has results in this direction.
A {\em binary structure} is one with binary vocabulary.

We say that a structure $\mcM$ is a {\em reduct} of a structure $\mcM'$ (possibly with another vocabulary)
if they have the same universe and for every positive integer $n$ and every relation $R \subseteq M^n$,
if $R$ is definable in $\mcM$ without parameters, then $R$ is definable in $\mcM'$ without parameters.
For any structure $\mcM$, $\mcM\meq$ denotes the extension of $\mcM$ by {\em imaginary elements} \cite{Hod, She}.
Note that understanding what kind of structures can be defined in $\mcM\meq$ is roughly the same as understanding 
which structures can be interpreted in $\mcM$.

We address the following problems:
{\em Suppose that $\mcM$ has finite relational vocabulary, is homogeneous and simple, $E \subset M$ is finite, $D \subseteq M\meq$ is
$E$-definable, only finitely many 1-types over $E$ are realized in $D$,  and for every $d \in D$ the SU-rank of the type of $d$ over $E$ is 1.}
\begin{itemize}
\item[(A)] What are the possible behaviours of the algebraic closure restricted to $D$ if elements of $E$ may be used as constants?

\item[(B)] Let $\mcD$ be the structure with universe $D$ which for every $n$ and $E$-definable $R \subseteq D^n$ has a relation symbol
which is interpreted as $R$ (and the vocabulary of $\mcD$ has no other symbols). We call $\mcD$ a {\em canonically embedded structure over $E$}. 
Note that the vocabulary of $\mcD$ is relational but {\em not} finite.
Now we ask whether $\mcD$ is necessarily a reduct of a homogeneous structure with finite relational vocabulary?
\end{itemize}

Macpherson \cite{Mac91} has shown that no infinite vector space over a finite field can be interpreted in a homogeneous structure
over a finite relational language,
which implies that, in~(A), the pregeometry of $D$ induced by the algebraic closure cannot be isomorphic to the pregeometry induced
by linear span in a vector space over a finite field.
If we assume, in addition to the assumptions made above (before~(A)), that $\mcM$ is one-based, then it follows from \cite{Mac91} and
work of De Piro and Kim \cite[Corollary~3.23]{PK} that algebraic closure restricted to $D$ is {\em trivial}, i.e.
if $d \in D$, $B \subseteq D$ and $d \in \acl(B \cup E)$, then there is $b \in B$ such that $d \in \acl(\{b\} \cup E)$.
But what if we do not assume that $\mcM$ is one-based?

Let $\mcD_0$ be the reduct of $\mcD$ to the relation symbols with arity at most 2.
If the answer to the question in~(B) is `yes' in the strong sense that $\mcD$ is a reduct of $\mcD_0$
and $\mcD_0$ is homogeneous, then 
Remark~\ref{remark on canonically embedded geometries}
below implies that algebraic closure and dependence restricted to $D$ are trivial.
If $\mcM$ is supersimple with finite SU-rank and the assumptions about $\mcD$ and $\mcD_0$ 
hold not only for this particular $\mcD$, but for
all $\mcD$, then it follows from \cite[Corollary~4.7]{HKP}, \cite[Corollary~3.23]{PK} and some additional straightforward arguments
that {\em the theory of  $\mcM$} has {\em trivial dependence} (Definition~\ref{definition of degenerate and trivial dependece} below).

In the other direction, if $\mcM$ is binary, has trivial dependence and $\acl(\{d\} \cup E) \cap D = \{d\}$ for
all $d \in D$ (so $D$ is a {\em geometry}), then,
by~Theorem~\ref{homogeneity of geometries},
$\mcD$ is a reduct of a binary homogeneous structure;
in fact $\mcD$ is a reduct of a {\em binary random structure} in the sense of Section~\ref{random structures}.

Thus we establish that, at least for binary $\mcM$, the problems~(A) and~(B) are closely related, 
although we do not know whether our partial conclusions to~(A) and~(B) in the binary case are equivalent.
Neither do we solve any one of problems~(A) or~(B). 
So in particular, the problem whether algebraic closure restricted to $D$ (and using constants from $E$) can be nontrivial for
some binary, homogeneous and simple $\mcM$ remains open.
Nevertheless, 
Theorem~\ref{homogeneity of geometries} is used in~\cite{Kop-one-based}
where a subclass of the countable, binary, homogeneous, simple and one-based structures is classified in a fairly concrete way;
namely the class of such structures which have {\em height} 1 in the sense of \cite{Djo06}, 
roughly meaning that the structure is ``coordinatized'' by a definable set of SU-rank 1.

This article is organized as follows.
In Section~\ref{Preliminaries} we recall definitions and results about homogeneous structures and simple structures,
in particular the independence theorem and consequences of $\omega $-categoricity and simplicity together, especially
with regard to imaginary elements. We also explain what is meant by a binary random structure.

In Section~\ref{Sets of rank one in simple homogeneous structures} 
we prove results implying that if $\mcM$ and $D$ are as assumed before~(A) above and $\mcM$ is binary,
then algebraic closure and dependence restricted to $D$ are trivial.
In Section~\ref{geometries} we prove the next main result,
Theorem~\ref{homogeneity of geometries}, saying that if  $\mcM$ and $D$ are as assumed before~(A), $\mcM$ is binary 
and its theory has trivial dependence, then $\mcD$ is a reduct of a binary random structure.
In order to prove Theorem~\ref{homogeneity of geometries} we use a more technical result, 
Theorem~\ref{homogeneity of D}, which is proved in Section~\ref{implications of trivial dependence}, 
where most of the technical 
(and simplicity theoretic) work is done.
The proofs assume a working knowledge in stability/simplicity theory, as can be found in \cite{Cas, Wag}.

\section{Preliminaries}\label{Preliminaries}

\subsection{General notation and terminology}

\noindent
A vocabulary (signature) is called {\em relational} if it only contains relation symbols.
For a finite relational vocabulary the maximal $k$ such that some relation symbol has arity $k$
is called its {\em maximal arity}.
If  $V$ is a {\em  finite} vocabulary and the maximal arity is 2 then we call $V$ {\em binary}
(although it may contain unary relation symbols), and in this case a $V$-structure is called a {\em binary structure}.
We denote (first-order) structures by $\mcA, \mcB, \ldots, \mcM, \mcN, \ldots$ and their respective universes by 
$A, B, \ldots, M, N, \ldots$. 
By the cardinality of a structure we mean the cardinality of its universe.
To emphasize the cardinality of a finite structure we sometimes call a structure with cardinality 
$k < \omega$ a {\em $k$-structure}, or {\em $k$-substructure} if it is seen as a substructure of some other structure.
Finite sequences (tuples) of elements of some structure (or set in general) will be denoted
$\bar{a}, \bar{b}, \ldots$, while $a, b, \ldots$ usually denote elements from the universe of some structure.
The notation $\bar{a} \in A$ means that every element in the sequence $\bar{a}$ belongs to $A$. 
Sometimes we write $\bar{a} \in A^n$ to show that the length of $\bar{a}$, denoted $|\bar{a}|$, is $n$ and all elements of  $\bar{a}$
belong to $A$.
By $\rng(\bar{a})$, the {\em range of $\bar{a}$}, we denote the set of elements that occur in $\bar{a}$.
In order to compress notation, we sometimes, in particular together with type notation and the symbol `$\ind$'
(for independence),  write `$AB$' instead of `$A \cup B$', or `$\bar{a}$' instead of `$\rng(\bar{a})$'.

Suppose that $\mcM$ is a structure, $A \subseteq M$ and $\bar{a} \in M$.
Then $\acl_\mcM(A)$, $\dcl_\mcM(A)$ and $\tp_\mcM(\bar{a} / A)$ denote the 
{\em algebraic closure} of $A$ with respect to $\mcM$, the {\em definable closure} of $A$ with respect to $\mcM$ and
the {\em complete type of $\bar{a}$ over $A$} with respect to $\mcM$, respectively
(see for example \cite{Hod} for definitions).
By $S_n^\mcM(A)$ we denote the set of all complete $n$-types over $A$ with respect to $\mcM$.
We abbreviate $\tp_\mcM(\bar{a} / \es)$ with $\tp_\mcM(\bar{a})$.
The notation $\acl_\mcM(\bar{a})$ is an abbreviation of $\acl_\mcM(\rng(\bar{a}))$, and similarly for `$\dcl$'.

We say that $\mcM$ is {\em $\omega$-categorical}, respectively {\em simple}, if $Th(\mcM)$ has that property,
where $Th(\mcM)$ is the complete theory of $\mcM$ (see \cite{Hod} and \cite{Cas, Wag} for definitions).
Let $A \subseteq M$ and $R \subseteq M^k$. We say that $R$ is {\em $A$-definable} (with respect to $\mcM$)
if there is a formula $\varphi(\bar{x}, \bar{y})$ (without parameters) and $\bar{a} \in A$ such that 
$R = \{\bar{b} \in M^k : \mcM \models \varphi(\bar{b}, \bar{a})\}$.
In this case we also denote $R$ by $\varphi(\mcM, \bar{a})$.
Similarly, for a type $p(\bar{x})$ (possibly with parameters)
we let $p(\mcM)$ be the set of all tuples of elements in $\mcM$
that realize $p$, and $\mcM \models p(\bar{a})$ means that $\bar{a}$ realizes $p$ in $\mcM$.
{\em Definable without parameters} means the same as $\es$-definable.

\begin{defin}\label{definition of reduct}{\rm
(i) If $\mcM$ is a structure with relational vocabulary and $A \subseteq M$, then $\mcM \uhrc A$ denotes the
substructure of $\mcM$ with universe $A$.\\
(ii) If $\mcM$ is a $V$-structure and $V' \subseteq V$, then $\mcM \uhrc V'$ denotes the {\em reduct of $\mcM$ to the vocabulary $V'$}.
}\end{defin}

\noindent
Note that if $\mcM$ is a $V$-structure and $V' \subseteq V$, then $\mcM \uhrc V'$ is a reduct of $\mcM$ in the sense defined in
Section~\ref{Introduction}.

\subsection{Homogeneity, Fra\"{i}ss\'{e} limits and $\omega$-categoricity}
\label{Homogeneous structures}

\begin{defin}\label{definition of homogeneous}{\rm
(i) Let $V$ be a relational vocabulary and $\mcM$ a $V$-structure. We call $\mcM$ {\em homogeneous} if its universe is countable and
for all finite substructures $\mcA$ and $\mcB$ of $\mcM$, every isomorphism from $\mcA$ to $\mcB$ can be extended to an
automorphism of $\mcM$.\\
(ii) A structure $\mcM$ (for any vocabulary) is called {\em $\omega$-homogeneous} if whenever $0 < n < \omega$,
$a_1, \ldots, a_n, a_{n+1}, b_1, \ldots, b_n \in M$ and $\tp(a_1, \ldots, a_n) = \tp(b_1, \ldots, b_n)$,
there is $b_{n+1} \in M$ such that $\tp(a_1, \ldots, a_{n+1}) = \tp(b_1, \ldots, b_{n+1})$.
}\end{defin}

\begin{defin}\label{definition of HP and AP}{\rm
Let $V$ be a finite relational vocabulary and let $\mbK$ be a class of finite $V$-structures which is {\em closed under isomorphism}, that is,
if $\mcA \in \mbK$ and $\mcB \cong \mcA$, then $\mcB \in \mbK$.\\
(i) $\mbK$ has the {\em hereditary property}, abbreviated HP, if $\mcA \subseteq \mcB \in \mbK$ implies that $\mcA \in \mbK$.\\
(ii) $\mbK$ has the {\em amalgamation property}, abbreviated AP, if the following holds:
if $\mcA, \mcB, \mcC \in \mbK$ and $f_\mcB : \mcA \to \mcB$ and $f_\mcC : \mcA \to \mcC$ are embeddings
then there are $\mcD \in \mbK$ and embeddings 
$g_\mcB : \mcB \to \mcD$ and $g_\mcC : \mcC \to \mcD$ such that $g_\mcB \circ f_\mcB = g_\mcC \circ f_\mcC$.\\ 
(iii) If $\mcM$ is a $V$-structure, then $\mathbf{Age}(\mcM)$ is the class of all $V$-structures that are isomorphic with some 
finite substructure of $\mcM$.
}\end{defin}

\noindent
We allow structures with empty universe if the vocabulary is relational (as generally assumed in this article), 
so if $\mbK$ has the hereditary property then the structure with empty universe belongs to $\mbK$.
It follows that if the vocabulary is relational then the {\em joint embedding property} \cite{Hod} is a consequence of the amalgamation property,
which is the reason why we need not bother about the former in the present context.
The following result of Fra\"{i}ss\'{e} (\cite{Fra54}, \cite{Hod} Theorems~7.1.2 and 7.1.7) 
relates homogeneous structures to finite structures, and shows how the 
former can be constructed from the later.

\begin{fact}\label{Fraisses theorem}
Let $V$ be a finite relational vocabulary. \\
(i) Suppose that $\mbK$ is a class of finite $V$-structures which is closed under
isomorphism and has HP and AP.
Then there is a unique, up to isomorphism, countable $V$-structure $\mcM$ such that
$\mcM$ is homogeneous and $\mathbf{Age}(\mcM) = \mbK$.\\
(ii) If $\mcM$ is a homogeneous $V$-structure, then $\mathbf{Age}(\mcM)$ has  HP and AP.
\end{fact}

\begin{defin}\label{Fraisse limit}{\rm
Suppose that $V$ is a finite relational vocabulary and that $\mbK$ is a class of finite $V$-structures
which is closed under isomorphism and has HP and AP. 
The unique (up to isomorphism) countable structure $\mcM$ such that $\mathbf{Age}(\mcM) = \mbK$ is 
called the {\em Fra\"{i}ss\'{e} limit} of $\mbK$.
}\end{defin}

\noindent
Part~(i) in the next fact is Corollary~7.4.2 in~\cite{Hod} (for example).
Part~(ii) follows from the well known characterisation of $\omega$-categorical structures by Engeler, Ryll-Nardzewski and Svenonius 
(\cite{Hod}, Theorem~7.3.1),
which will frequently be used without further reference.
Part~(iii) follows from a straightforward back and forth argument.

\begin{fact}\label{equivalences to being homogeneous}
Let $V$ be a relational vocabulary and $\mcM$ an infinite countable $V$-structure.\\
(i) If $V$ is finite then $\mcM$ is homogeneous if and only if $\mcM$ is $\omega$-categorical and has elimination of quantifiers.\\
(ii) If $\mcM$ is $\omega$-categorical, then $\mcM$ is $\omega$-saturated and $\omega$-homogeneous.\\
(iii) If $\mcM$ is countable and $\omega$-homogeneous, then the following holds:
if $0 < n < \omega$, $a_1, \ldots, a_n, b_1, \ldots, b_n \in M$ and 
$\tp_\mcM(a_1, \ldots, a_n) = \tp_\mcM(b_1, \ldots, b_n)$, then there is an automorphism $f$ of $\mcM$ 
such that $f(a_i) = b_i$ for all $i$.
\end{fact}

\subsection{Binary random structures}\label{random structures}

Let $V$ be a binary vocabulary (and therefore finite).

\begin{defin}\label{definition of 1-adequate}{\rm
A class $\mbK$ of finite $V$-structures is called {\em 1-adequate} if it has HP and the following property with respect to 
1-structures:
\begin{itemize}
\item[] If $\mcA, \mcB \in \mbK$ are 1-structures, then there is $\mcC \in \mbK$ such that $\mcA \subseteq \mcC$ and $\mcB \subseteq \mcC$.
\end{itemize}
}\end{defin}

\noindent
{\bf Construction of a binary random structure:}
Let $\mbP_2$ be a 1-adequate class of $V$-structures such that  $\mbP_2$ contains a 2-structure.
We think of $\mbP_2$ as containing the isomorphism types of ``permitted'' 1-(sub)structures and 2-(sub)structures.
Then let $\mbR\mbP_2$ be the class of all finite $V$-structures $\mcA$ such that for $k = 1,2$ 
every $k$-substructure of $\mcA$ is isomorphic to some
member of $\mbP_2$. 
Obviously, $\mbR\mbP_2$ has HP, because the 1-adequateness of $\mbP_2$ implies that $\mbP_2$ has HP.
The 1-adequateness of $\mbP_2$  implies that any two 1-structures of $\mbR\mbP_2$ can be embedded into
a 2-structure of $\mbP_2$. From this it easily follows that $\mbR\mbP_2$ has AP.
Let $\mcF$ be the Fra\"{i}ss\'{e} limit of $\mbR\mbP_2$.
We call $\mcF$ the {\em random structure over $\mbP_2$}, or more generally a {\em binary random structure}.
This is motivated by the remark below. 
But first we show that the well known ``random graph'' (or ``Rado graph'')  is a binary random structure in this sense.

\begin{exam}\label{example of the random graph}{\rm ({\bf The random graph})
Let $V = \{R\}$, where $R$ is a binary relation symbol and let $\mbP_2$ be the following class (in fact a set) of
$V$-structures, where $(A, B)$ denotes the $\{R\}$-structure with universe $A$ and where $R$ is intepreted as $B \subseteq A^2$:
\[ \mbP_2 \ = \ \big\{ (\es, \es),  (\{1\}, \es), (\{1,2\}, \es), (\{1,2\}, \{(1,2), (2,1)\}) \big\} \]
If $\mbR\mbP_2$ is as in the construction above,
then $\mbR\mbP_2$ is the class of all finite undirected graphs (without loops), which has HP and AP, and the 
Fra\"{i}ss\'{e} limit of it is (in a model theoretic context) often called the {\em random graph}.
}\end{exam}

\begin{rem}\label{remark about random structure over P_2}{\rm {\bf (Random structures and zero-one laws)}
Let $\mbP_2$ and $\mbR\mbP_2$ be as in the construction of a binary random structure above.
Then $\mbR\mbP_2$ is a {\em parametric class} in the sense of Definition~4.2.1 in \cite{EF} or Section~2 of \cite{Ober}.
Hence, by Theorem~4.2.3 in \cite{EF}, $\mbR\mbP_2$ has a (labelled) 0-1 law (with the uniform probability measure).
This is proved  by showing that all extension axioms that are compatible with $\mbR\mbP_2$ hold with probability
approaching 1 as the (finite) cardinality of members of $\mbR\mbP_2$ approaches infinity; see statement~(5) on page~76 in \cite{EF}.
(Alternatively, one can use the terminology of \cite{Kop12} and show that $\mbR\mbP_2$ ``admits $k$-substitutions'' for every positive
integer $k$, and then apply Theorem~3.15 in \cite{Kop12}.)
It follows that if $T_{\mbR\mbP_2}$ is the set of all $V$-sentences $\varphi$ with asymptotic probability~1 (in $\mbR\mbP_2$), 
then all extension axioms that are compatible with $\mbR\mbP_2$ belong to $T_{\mbR\mbP_2}$.
Let $\mcF$ be the Fra\"{i}ss\'{e} limit of $\mbR\mbP_2$.
Then $\mcF$ satisfies every extension axiom which is compatible with $\mbR\mbP_2$
(since if $\mcA \subseteq \mcF$ and $\mcA \subseteq \mcB \in \mbR\mbP_2$, then there is an embedding of
$\mcB$ into $\mcF$ which is the identity on $\mcA$).
By a standard back-and-forth argument, it follows that if $\mcM$ is a countable model of $T_{\mbR\mbP_2}$, 
then $\mcM \cong \mcF$ and hence $\mcF \models T_{\mbR\mbP_2}$.
}\end{rem}

\noindent
The construction of a binary random structure can of course be generalised to any finite relational 
(not necessarity binary) vocabulary.

\subsection{Simple $\omega$-categorical structures, imaginary elements and rank}\label{Simple omega-categorical structures}

\noindent
We will work with concepts from stability/simplicty theory, including imaginary elements. 
That is, we work in the
structure $\mcM\meq$ obtained from a structure $\mcM$ by adding ``imaginary'' elements, in the way explained in \cite{Hod, She}, for example.
In the case of $\omega$-categorical simple theories, some notions and results become easier than in the
general case. For example, every $\omega$-categorical simple theory has elimination of
hyperimaginaries, so we need not consider ``hyperimaginary elements'' or the ``bounded closure'';
it suffices to consider imaginary elements and algebraic closure, so we need not go beyond $\mcM\meq$.
The results about $\omega$-categorical simple structures that will be used, 
often without explicit reference, are stated below, with proofs or at least indications of how they follow from well known results
in stability/simplicity theory or model theory in general.

Let $\mcM$ be a $V$-structure. 
Although we assume familiarity with  $\mcM\meq$, the universe of which is denoted $M\meq$,  
we recall part of its construction (as in \cite{Hod, She} for instance), 
since the distinction between different ``sorts'' of elements of $\mcM\meq$ matters in the present work.
For every $0 < n < \omega$ and every equivalence relation $E$ on $M^n$ which is $\es$-definable in $\mcM$,
$V\meq$ (the vocabulary of $\mcM\meq$) contains a unary relation symbol $P_E$ and an $(n+1)$-ary relation
symbol $F_E$ (which do not belong to $V$). $P_E$ is
interpreted as the set of $E$-equivalence classes and, for all $\bar{a} \in (M\meq)^n$ and each $c \in M\meq$,
$\mcM\meq \models F_E(\bar{a}, c)$ if and only if $\bar{a} \in M^n$, $c$ is an $E$-equivalence class and $\bar{a}$ belongs to $c$. 
(So the interpretation of $F_E$ is the graph of a function from $M^n$ to the set of all $E$-equivalence classes.)
{\em The notation $F(\bar{a}, c)$ means that $\mcM\meq \models F_E(\bar{a}, c)$ for some $n$ and some $\es$-definable
equivalence relation $E$ on $M^n$.}

A {\em sort} of $\mcM\meq$ is, by definition, a set of the form $S_E = \{a \in M\meq : \mcM\meq \models P_E(a) \}$ for some 
$E$ as above. If $A \subseteq M\meq$ and there are only finitely many $E$ such that $A \cap S_E \neq \es$, 
then we say that {\em only finitely many sorts are represented in $A$}.
Note that `=', the identity relation, is an $\es$-definable equivalence relation on $M$ and every $=$-class is a singleton.
Therefore $M$ can (and will) be identified with the sort $S_=$, which we call the {\em real sort}.
Hence every element of $M\meq$ belongs to $S_E$ for some $E$.
If $\mcN \equiv \mcM\meq$ then every element $a \in N$ such that $\mcN \models P_=(a)$ is 
called a {\em real element} and every element $a \in N$ such that $\mcN \models P_E$ for some $E$
is called an {\em imaginary element} (so real elements are special cases of imaginary elements). 
However, the set 
\[ \{\neg P_E(x) : \text{ $E$ is a $\es$-definable equivalence relation on $M^n$ for some $n$} \} \]
is consistent with $Th(\mcM\meq)$ (by compactness), so some model of $Th(\mcM\meq)$ will contain elements which are neither real nor imaginary. 
This also shows that $\mcM\meq$ is {\em not} $\omega$-saturated even if 
$\mcM$ is (which is the case if $\mcM$ is $\omega$-categorical).
However, if $\mcM$ is $\omega$-categorical and $A \subseteq M\meq$ is finite, then every type $p \in S_n^{\mcM\meq}(A)$
which is realized by an  $n$-tuple of imaginary elements in some elementary extension of $\mcM\meq$ is already realized in $\mcM\meq$, 
as stated by Fact~\ref{second fact about canonically embedded structures} below.
The first fact below follows from Theorem~4.3.3 in \cite{Hod} or Lemma~III.6.4 in \cite{She}.

\begin{fact}\label{Meq does not change truth}
For all $\bar{a}, \bar{b} \in M$, 
$\tp_\mcM(\bar{a}) = \tp_\mcM(\bar{b})$ if and only if $\tp_{\mcM\meq}(\bar{a}) = \tp_{\mcM\meq}(\bar{b})$.
\end{fact}

\begin{fact}\label{only finitely many realizations}
Suppose that $\mcM$ is $\omega$-categorical,
let $A \subseteq M\meq$ and suppose that only finitely many sorts are represented in $A$. \\
(i) For every $n < \omega$ and finite $B \subseteq \mcM\meq$, only finitely many types from $S_n^{\mcM\meq}(\acl_{\mcM\meq}(B))$
are realized by $n$-tuples in $A^n$.\\
(ii) For every $n < \omega$ and finite $B \subseteq \mcM\meq$, $\acl_{\mcM\meq}(B) \cap A$ is finite.
\end{fact}

\noindent
{\em Proof.}
Let $B' \subseteq M$ be finite and such that $B \subseteq \acl_{\mcM\meq}(B')$.
By $\omega$-categoricity, there are, up to equivalence in $\mcM$, only finitely many formulas in free variables $x_1, \ldots, x_n$
with parameters from $B'$, so part~(i) is a consequence of Lemma~6.4 of Chapter~III in \cite{She} (or use 
Theorem~4.3.3 in \cite{Hod}).
Part~(ii) follows from part~(i).
\hfill $\square$

\begin{defin}\label{definition of canonically embedded structure}{\rm
Suppose that $A \subseteq M\meq$ is finite.
We say that a structure $\mcN$ is {\em canonically embedded in $\mcM\meq$ over $A$} if $N$ is an $A$-definable subset of $M\meq$ and
for every $0 < n < \omega$ and every relation $R \subseteq N^n$ which is $A$-definable in $\mcM\meq$ there is a
relation symbol in the vocabulary of $\mcN$ which is interpreted as $R$ and the vocabulary of $\mcN$ contains no other relation symbols
(and no constant or function symbols).
}\end{defin}

\noindent
The following is immediate from the definition:

\begin{fact}\label{first fact about canonically embedded structures}
If $A \subseteq M\meq$ is finite and $\mcN$ is canonically embedded in $\mcM\meq$ over $A$, then for all $\bar{a}, \bar{b} \in N$ and all $C \subseteq N$,
$\acl_\mcN(C) = \acl_{\mcM\meq}(CA) \cap N$ and
$\tp_\mcN(\bar{a} / C) = \tp_\mcN(\bar{b} / C)$ if and only if $\tp_{\mcM\meq}(\bar{a} / CA) = \tp_{\mcM\meq}(\bar{b} / CA)$. 
\end{fact}

\begin{fact}\label{second fact about canonically embedded structures}
Suppose that $\mcM$ is $\omega$-categorical. \\
(i) If $\mcN$ is canonically embedded in $\mcM\meq$ over a finite $A \subseteq M\meq$ and only finitely sorts are 
represented in $N$, then $\mcN$ is $\omega$-categorical and therefore $\omega$-saturated.\\
(ii) If $A \subseteq M\meq$ is finite and $\bar{a} \in M\meq$,
then $\tp_{\mcM\meq}(\bar{a} / \acl_{\mcM\meq}(A))$ is isolated. \\
(iii) If $A \subseteq \mcM\meq$ is finite, $n < \omega$ and $p \in S_n^{\mcM\meq}(\acl_{\mcM\meq}(A))$ is realized 
in some elementary extension of $\mcM\meq$ by an $n$-tuple of imaginary elements, 
then $p$ is realized in $\mcM\meq$.\\
(iv) If $\mcM$ is countable, then $\mcM\meq$ is $\omega$-homogeneous.
\end{fact}

\noindent
{\em Proof.}
(i) If $\mcM$ is $\omega$-categorical, then, by the characterization of its complete theory by Engeler, Ryll-Nardzewski and Svenonius
(the characterisation by isolated types),
Fact~\ref{first fact about canonically embedded structures} and, for example,  Lemma~6.4 of Chapter~III in \cite{She} (or Fact 1.1 in \cite{Djo06}),
it follows that $\mcN$ is $\omega$-categorical (and hence $\omega$-saturated).

(ii) For $\omega$-categorical $\mcM$, finite $A \subseteq M\meq$ and $\bar{a} \in M\meq$, 
it follows from Fact~\ref{first fact about canonically embedded structures} 
and part~(i) that $\tp(\bar{a} / A)$ is isolated. From the assumption that $\tp(\bar{a} / \acl_{\mcM\meq}(A))$
is not isolated it is straightforward to derive a contradiction to Fact~\ref{only finitely many realizations}.
Parts~(iii) and~(iv) follow from part~(ii).
\hfill $\square$
\\

\noindent
Besides the above stated consequences of $\omega$-categoricity, the proofs in
Sections~\ref{Sets of rank one in simple homogeneous structures}
and~\ref{implications of trivial dependence}
use the so called {\em independence theorem} for simple theories \cite{Cas, Wag}.
Every $\omega$-categorical simple theory has elimination of hyperimaginaries and, with respect to it, `Lascar strong types' are equivalent
with strong types (\cite{Cas} Theorem~18.14, \cite{Wag} Lemma~6.1.11), from which it follows that any two finite tuples 
$\bar{a}, \bar{b} \in M\meq$ have the same Lascar strong type over
a finite set $A \subseteq M\meq$ if and only if they have the same type over $\acl_{\mcM\meq}(A)$.
Therefore the independence theorem implies the following, which is the version of it that we will use:

\begin{fact}\label{independence theorem} {\bf (The independence theorem for simple $\omega$-categorical structures and finite sets)}
Let $\mcM$ be a simple and $\omega$-categorical structure and
let $A, B, C \subseteq \mcM\meq$ be finite.
Suppose that $B \underset{A}{\ind} C$,
$n < \omega$, $\bar{b}, \bar{c} \in (M\meq)^n$, 
$\bar{b} \underset{A}{\ind} B$, $\bar{c} \underset{A}{\ind} C$
and
\[ \tp_{\mcM\meq}(\bar{b} / \acl_{\mcM\meq}(A)) \ = \ 
\tp_{\mcM\meq}(\bar{c} / \acl_{\mcM\meq}(A)). \] 
Then there is $\bar{d} \in (M\meq)^n$ such that 
\begin{align*} 
&\text{$\tp_{\mcM\meq}(\bar{d} / B \cup \acl_{\mcM\meq}(A)) \ = \ \tp_{\mcM\meq}(\bar{b} / B \cup \acl_{\mcM\meq}(A))$,} \\
&\text{$\tp_{\mcM\meq}(\bar{d} / C \cup \acl_{\mcM\meq}(A)) \ = \ \tp_{\mcM\meq}(\bar{c} / C \cup \acl_{\mcM\meq}(A))$} 
\end{align*}
and $\bar{d}$ is independent from $B \cup C$ over $A$.
\end{fact}

\noindent
By induction one easily gets the following, which is sometimes more practical:

\begin{cor}\label{corollary to independence theorem}
Let $\mcM$ be a simple and $\omega$-categorical structure, $2 \leq k < \omega$ and
let $A, B_1, \ldots, B_k \subseteq \mcM\meq$ be finite.
Suppose that $\{B_1, \ldots, B_k\}$ is independent over $A$,
$n < \omega$, $\bar{b}_1, \ldots, \bar{b}_k \in (M\meq)^n$
and, for all $i, j \in \{1, \ldots, k\}$, $\bar{b}_i \underset{A}{\ind} B_i$ and
\[ \tp_{\mcM\meq}(\bar{b}_i / \acl_{\mcM\meq}(A)) \ = \ \tp_{\mcM\meq}(\bar{b}_j / \acl_{\mcM\meq}(A)). \]
Then there is $\bar{b} \in (M\meq)^n$ such that, for all $i = 1, \ldots, k$,
\[ \tp_{\mcM\meq}(\bar{b} / B_i \cup \acl_{\mcM\meq}(A)) \ = \ \tp_{\mcM\meq}(\bar{b}_i / B_i \cup \acl_{\mcM\meq}(A)) \]
and $\bar{b}$ is independent from $B_1 \cup \ldots \cup B_k$ over $A$.
\end{cor}

\noindent
Suppose that $T$ is a simple theory. 
For every type $p$ (possibly over a set of parameters) with respect to $T$, there is a notion of {\em $\su$-rank} of $p$, denoted $\su(p)$;
it is defined in \cite{Cas, Wag} for instance. We abbreviate $\su(\tp_\mcM(\bar{a} / A))$ with $\su(\bar{a} / A)$.
For any type $p$,  $\su(p)$ is either ordinal valued or undefined (or alternatively given the value $\infty$).

\section{Sets of rank one in simple homogeneous structures}\label{Sets of rank one in simple homogeneous structures}

\noindent
In this section we derive consequences for sets with rank one in simple homogeneous structures with 
the $n$-dimensional amalgamation property for strong types (defined below), where $n$ is the maximal arity of the vocabulary.
A consequence of the independence theorem is that all simple structures have the 2-dimensional amalgamation property for strong types.
We will use the notation $\mcP(S)$ for the powerset of the set $S$, and let $\mcP^-(S) = \mcP(S) \setminus \{S\}$. 
Every $n < \omega$ is identified with the set $\{0, \ldots, n-1\}$, and hence the notation $\mcP(n)$ makes sense. 
For a type $p$, $\dom(p)$ denotes the set of all parameters that occur in formulas in $p$.
We now consider the `strong $n$-dimensional amalgamation property for Lascar strong types', 
studied by Kolesnikov in \cite{Kol05} (Definition~4.5). 
However, we only need it for real elements 
and in the present context `Lascar strong type' is the same as `type over an algebraically closed set'.

\begin{defin}{\rm
Let $T$ be an $\omega$-categorical  and simple theory and let $n < \omega$. \\
(i) A set of types $\{p_s(\bar x) | s\in \mcP^-(n)\}$ (with respect to $\mcM\meq$ for some $\mcM \models T$) 
is called an {\em $n$-independent system of strong types over $A$} (where $A \subseteq \mcM\meq$)
if it satisfies the following properties:
\begin{itemize}
\item $\dom(p_\emptyset) = A$.
\item for all $s,t \in \mcP^-(n)$ such that $s\subseteq t$, $p_t$ is a nondividing extension of $p_s$.
\item for all $s,t \in \mcP^-(n)$, $\dom(p_s)\us{\dom(p_{s\cap t})}{\ind} \dom(p_t)$.
\item for all $s,t \in \mcP^-(n)$, $p_s$ and $p_t$ extend the same type over $\acl_{\mcM\meq}(\dom(p_{s\cap t}))$.
\end{itemize}
(ii) We say that $T$ (and any $\mcN \models T$) has the {\em  n-dimensional amalgamation property for strong types} 
if for every $\mcM \models T$ and every
$n$-independent system of strong types $\{p_s(\bar x) | s\in \mcP^-(n)\}$ over some set $A \subseteq M\meq$,
there is a type $p^*$ such that $p^*$ extends $p_s$ for each $s \in \mcP^-(n)$ and
$p^*$ does not divide over $\bigcup_{s \in \mcP^-(n)}\dom(p_s)$.
}\end{defin}

\begin{rem}\label{remark that simplicity implies the 2-dimensional amalgamation property}{\rm
By the independence theorem (in the general case when the sets of parameters of the given types may be infinite \cite{Cas, Wag}), 
every $\omega$-categorical and simple theory has the 2-dimensional amalgamation property for strong types.
}\end{rem}

\begin{theor}\label{pregeometries are degenerate if the structure is homogeneous}
Suppose that $\mcM$ has a finite relational vocabulary with maximal arity $\rho$.
Also assume that $\mcM$ is countable, homogeneous and simple and has the $\rho$-dimensional amalgamation property for strong types.
Let $D, E \subseteq M$ where $E$ is finite, $D$ is $E$-definable, and  $\su(a / E) = 1$ for every $a \in D$.
If $a \in D$, $B \subseteq D$ and $a \in \acl_{\mcM\meq}(BE)$, then $a \in \acl_{\mcM\meq}(B' E)$ for some $B' \subseteq B$ with $|B'| < \rho$.
\end{theor}

\noindent
{\em Proof.}
Assume that $a \in D$, $B \subseteq D$ and $a \in \acl_{\mcM\meq}(BE)$.
Without loss of generality we may assume that $B$ is finite.
By induction on $|B|$ we prove that there is $B' \subseteq B$ such that $|B'| < \rho$ and $a \in \acl_{\mcM\meq}(B' E)$.
The base case is when $|B| < \rho$ and we are automatically done.

So suppose that $|B| \geq \rho$.
If $B$ is not independent over $E$ then there is $b \in B$ such that
$b \underset{E}{\nind} (B \setminus \{b\})$ and as $\su(b / E) = 1$ (by assumption)
we get $b \in \acl_{\mcM\meq}\big((B \setminus \{b\}) \cup E\big)$.
Hence $B \subseteq \acl_{\mcM\meq}(B'E)$ where $B' = B \setminus \{b\}$ is a proper subset of $B$,
so by the induction hypothesis we are done.

So now suppose, in addition, that $B$ is independent over $E$.
If $a \in \acl_{\mcM\meq}(B' E)$ for some proper subset $B' \subset B$, then we are done by the induction hypothesis.
Therefore assume, in addition, that $a \notin \acl_{\mcM\meq}(B' E)$ for every proper subset $B' \subset B$.

Let $n = |B|$, so $n \geq \rho$ and enumerate $B$ as $B = \{b_0, \ldots, b_{n-1}\}$. 
For each $S \in \mcP^-(\rho)$, let 
\[B_S = \acl_{\mcM\meq}\big(\{b_t : t \in S\} \cup  \{b_\rho, \ldots, b_{n-1}\} \cup E\big).\]
From the assumptions that $B$ is independent over $E$ and $a \notin \acl_{\mcM\meq}(B' E)$ for every proper subset $B' \subset B$
it follows that the types $\tp(a / B_S)$ form a $\rho$-independent system of strong types over 
$\acl_{\mcM\meq}(E \cup \{b_\rho, \ldots, b_{n-1}\})$.
As $Th(\mcM)$ has the $\rho$-dimensional amalgamation property for strong types 
(and using Fact~\ref{second fact about canonically embedded structures}), 
there is $a' \in D$ such that for every $i \in \{0, \ldots, \rho - 1\}$ and $S_i = \{0, \ldots, \rho - 1\} \setminus \{i\}$,
\begin{align}\label{a' satisfying amalgamation}
&\tp_{\mcM\meq}(a' / B_{S_i}) \ = \ \tp_{\mcM\meq}(a / B_{S_i}) \ \text{ and} \\
&a' \notin \acl_{\mcM\meq}(BE). \nonumber
\end{align}

\bigskip
\noindent
{\bf Claim.} The bijection $f : \mcM \uhrc BEa' \to \mcM \uhrc BEa$ defined by $f(b) = b$ for all $b \in BE$ 
and $f(a') = a$ is an isomorphism.

\medskip

\noindent
{\bf Proof of the Claim.}
Let $R$ be a relation symbol of the vocabulary of $\mcM$, so the arity of $R$ is at most $\rho$.
It suffices to show that if $\bar{a} \in BEa'$ then $\mcM \models R(\bar{a})$ 
if and only if $\mcM \models R(f(\bar{a}))$.
But this is immediate from~(\ref{a' satisfying amalgamation}) and the definition of $f$.
\hfill $\square$

\bigskip

\noindent
Since $\mcM$ is homogeneous and $B$ and $E$ are finite, there is an automorphism $g$ of $\mcM$ which extends $f$ from the claim.
Then  $g(a') = a$ and $g$ fixes $BE$ pointwise.
But since $a \in \acl_{\mcM\meq}(BE)$ and, by~(\ref{a' satisfying amalgamation}), $a' \notin \acl_{\mcM\meq}(BE)$ this contradicts that $g$ is an automorphism.
\hfill $\square$
\\

\noindent
By using the previous theorem and Remark~\ref{remark that simplicity implies the 2-dimensional amalgamation property}
we get the following:

\begin{cor}\label{corollary giving trivial algebraic closure if M is binary}
Suppose that $\mcM$ is a countable, binary, homogeneous and simple structure.
Let $D, E \subseteq M$ where $E$ is finite, $D$ is $E$-definable and $\su(a / E) = 1$ for every $a \in D$.
If $a \in D$, $B \subseteq D$ and $a \in \acl_{\mcM\meq}(BE)$, then $a \in \acl_{\mcM\meq}(\{b\} \cup E)$ for some $b \in B$.
\end{cor}

\begin{defin}\label{definition of degenerate and trivial dependece}{\rm
Let $T$ be a simple theory. \\
(i) Suppose that $\mcM \models T$ and $E \subseteq M$.
We say that $D \subseteq M\meq$ has  {\em $n$-degenerate dependence over $E$} 
if for all $A,B,C\subseteq D$ such that $A\us{CE}{\nind} B$ there
is $B_0\subseteq B$ such that $|B_0|\leq n$ and $A\us{CE}{\nind} B_0$. \\
(ii) We say that $T$ has {\em trivial dependence} if  whenever $\mcM \models T$,
$A, B, C_1, C_2 \subseteq M\meq$ and $A \underset{B}{\nind}C_1 C_2$, then $A \underset{B}{\nind} C_i$ for $i = 1$ or $i = 2$.
A simple structure $\mcM$ has {\em trivial dependence} if its complete theory $Th(\mcM)$ has it.
}\end{defin}

\begin{theor}\label{dependence is rho-degenerated if M has the rho-dimensional amalgamation property}
Suppose that $\mcM$ has a finite relational vocabulary with maximal arity~$\rho$.
Also assume that $\mcM$ is countable, homogeneous and simple and has the $\rho$-dimensional amalgamation property for strong types.
Let $D, E \subseteq M$ where $E$ is finite,  $D$ is $E$-definable and $\su(a / E) = 1$ for every $a \in D$.
Then $D$ has $(\rho - 1)$-degenerate dependence over $E$.
\end{theor}

\noindent
{\em Proof.}
This is essentially an application of Theorem~\ref{pregeometries are degenerate if the structure is homogeneous},
basic properties of SU-rank and the Lascar (in)equalities (see for example Chapter 5.1 in \cite{Wag}, in particular Theorem~5.1.6).

Suppose that $B, C \subseteq D$, $\bar{a} = (a_1, \ldots, a_n) \in D^n$ and 
\begin{equation}\label{bar-a is depent of B over CE}
\bar{a} \underset{CE}{\nind} B. 
\end{equation}
If $\bar{a}$ is not independent over $CE$, then
(since $\su(d/E) = 1$ for all $d \in D$) there is a proper subsequence $\bar{a}'$ of $\bar{a}$ such that
$\rng(\bar{a}) \subseteq \acl(\bar{a}' CE)$ and hence $\bar{a}' \underset{CE}{\nind} B$.
If, in addition, $B' \subsetneq B$ and $\bar{a}' \underset{CE}{\nind} B'$, then $\bar{a} \underset{CE}{\nind} B'$.
Therefore we may assume that
\begin{equation}\label{bar-a is independet over CE}
\text{$\bar{a}$ is independent over $CE$.}
\end{equation}
Moreover, we may assume that 
\begin{equation}\label{a-i has SU-rank 1 over CE}
\su(a_i / CE) = 1 \ \ \text{ for every } i.
\end{equation}
For otherwise, $a_i \in \acl(CE)$ for some $i$, which implies that $\bar{a}$ is not independent 
over $CE$, contradicting~(\ref{bar-a is independet over CE}).

Now~(\ref{bar-a is independet over CE}),~(\ref{a-i has SU-rank 1 over CE}) and the Lascar equalities (for finite ranks) give
\begin{equation}\label{SU-rank of bar-a over CE is n}
\su(\bar{a} / CE) = n.
\end{equation}
Then~(\ref{bar-a is depent of B over CE}) and~(\ref{SU-rank of bar-a over CE is n}) (together with Lemma~5.1.4 in \cite{Wag} for example)
give
\[
\su(\bar{a} / BCE) < n,
\]
so $\bar{a}$ is not independent over $BCE$ and hence there is $i$ such that
\[
a_i \underset{BCE}{\nind} (\{a_1, \ldots, a_n\} \setminus \{a_i\}),
\]
so (by monotonicity of dependence)
\[
a_i \underset{CE}{\nind} B \cup (\{a_1, \ldots, a_n\} \setminus \{a_i\}),
\]
and by~(\ref{a-i has SU-rank 1 over CE}),
\[
a_i \in \acl(BCE \cup (\{a_1, \ldots, a_n\} \setminus \{a_i\})).
\]
By Theorem~\ref{pregeometries are degenerate if the structure is homogeneous},
there is $X \subseteq BC \cup (\{a_1, \ldots, a_n\} \setminus \{a_i\})$ such that
$|X| < \rho$ and $a_i \in \acl(XE)$.
Then
\[a_i \in \acl\big((X \cap B) \cup CE \cup (\{a_1, \ldots, a_n\} \setminus \{a_i\})\big).\]
Let $\bar{a}'$ be the proper subsequence of $\bar{a}$ in which $a_i$ removed.
Then, using~(\ref{SU-rank of bar-a over CE is n}), 
\begin{equation*}
\su(\bar{a} / (X \cap B) \cup CE) \ = \ \su(\bar{a}' / (X \cap B) \cup CE) \ < \ n \ =
\ \su(\bar{a} /CE),
\end{equation*}and hence
$\bar{a} \underset{CE}{\nind} (X \cap B)$ where $|X \cap B| < \rho$.
\hfill $\square$

\begin{rem}\label{remark about adding constants}{\rm
Suppose that $\mcM$ is homogeneous and simple and that $E \subseteq M$ is finite.
Let $\mcM_E$ be the expansion of $\mcM$ with a unary relation symbol $P_e$ for each $e \in E$ and interpret
$P_e$ as $\{e\}$. It is straightforward to verify that any isomorphism between finite substructures of $\mcM_E$
can be extended to an automorphism of $\mcM_E$, so it is homogeneous.
Moreover, since the notion of simplicity only depends on which relations are definable with parameters and exactly 
the same relations are definable with parameters in $\mcM_E$ as in $\mcM$ it follows that $\mcM_E$ is simple
(see for example \cite[Remark 2.26]{Cas}).
For the same reason, if $\mcM$ has trivial dependence, then so has $\mcM_E$.
}\end{rem}

\begin{cor}\label{corollary on canonically embedded geometries}
Suppose that $\mcM$ is countable, homogeneous, simple with a finite relational vocabulary with maximal arity $\rho$,
and with the $\rho$-dimensional amalgamation property for strong types.
Let $D \subseteq M\meq$ be $E$-definable for finite $E \subseteq M$, suppose that only finitely many sorts are represented in $D$
and that $\su(d / E) = 1$ for all $d \in D$.
Moreover, suppose that
if $n < \omega$, $a_1, \ldots, a_n, b_1, \ldots, b_n \in D$ and 
\begin{align*}
&\tp_{\mcM\meq}(a_{i_1}, \ldots, a_{i_\rho} /E) \ = \ \tp_{\mcM\meq}(b_{i_1}, \ldots, b_{i_\rho}/E) \ \ 
\text{ for all $i_1, \ldots, i_\rho \in \{1, \ldots, n\}$,} \\
&\text{then $\tp_{\mcM\meq}(a_1, \ldots, a_n/E) \ = \ \tp_{\mcM\meq}(b_1, \ldots, b_n/E)$.}
\end{align*}
Then $D$ has $(\rho - 1)$-degenerate dependence over $E$. 
\end{cor}

\noindent
{\em Proof.}
Suppose that $\mcM$, $D$ and $E$ satisfy the assumptions of the  corollary.
Let $\mcM_E$ be the expansion of $\mcM$ by a unary relation symbol $P_e$ for every $e \in E$ where $P_e$ is interpreted as $\{e\}$. 
By Remark~\ref{remark about adding constants}, $\mcM_E$ is homogeneous and simple.
Moreover, $D$ is $\es$-definable in $(\mcM_E)\meq$, so it is the universe of a canonically embedded structure $\mcD$ in $(\mcM_E)\meq$ over $\es$.
By~\ref{second fact about canonically embedded structures}, $\mcD$ is $\omega$-categorical and hence $\omega$-homogeneous.
As $\mcM$, and thus $\mcD$, is countable it follows that if $0 < n < \omega$, 
$a_1, \ldots, a_n, b_1, \ldots, b_n \in D$ and
$\tp_{\mcD}(a_1, \ldots, a_n) = \tp_{\mcD}(b_1, \ldots, b_n)$,
then there is an automorphism $f$ of $\mcD$ such that $f(a_i) = b_i$ for all $i$.
By assumption and Fact~\ref{Meq does not change truth}, if $n < \omega$, $a_1, \ldots, a_n, b_1, \ldots, b_n \in D$ and 
$\tp_{\mcD}(a_{i_1}, \ldots, a_{i_\rho}) = \tp_{\mcD}(b_{i_1}, \ldots, b_{i_\rho})$
for all $i_1, \ldots, i_\rho \in \{1, \ldots, n\}$, then 
\[\tp_{\mcD}(a_1, \ldots, a_n) = \tp_{\mcD}(b_1, \ldots, b_n).\]
Hence the reduct $\mcD_0$ of $\mcD$ to the relation symbols of arity at most $\rho$ is homogeneous (and $\mcD$ is a reduct of $\mcD_0$).
Clearly, $D$ is an $\es$-definable subset in $\mcD$ and by Fact~\ref{first fact about canonically embedded structures} we have
$\su(d / \es) = 1$ for all $d \in D$ when `$\su$' is computed in $\mcD$, as well as in $\mcD_0$ (since $\mcD$ is a reduct of $\mcD_0$).
Hence, Theorem~\ref{dependence is rho-degenerated if M has the rho-dimensional amalgamation property}
applied to $\mcM = \mcD_0$ implies that $D$
has $(\rho - 1)$-degenerate dependence when we consider $D$ as a $\es$-definable set within $\mcD_0$, and hence when $D$ is
considered as a $\es$-definable set within $\mcD$. 
From Fact~\ref{first fact about canonically embedded structures} 
it follows that $D$ has $(\rho - 1)$-degenerate dependence when we consider $D$ as a  $\es$-definable set within $(\mcM_E)\meq$.
It follows that
$D$ has $(\rho - 1)$-definable dependence over $E$ when we consider $D$ as an $E$-definable set within $\mcM\meq$.
\hfill $\square$

\begin{rem}\label{remark on canonically embedded geometries}{\rm
As mentioned earlier, every $\omega$-categorical simple theory has the $2$-dimen\-sional amalgamation property
for strong types.
So if $\mcM$ in Corollary~\ref{corollary on canonically embedded geometries} is binary, that is, if $\rho \leq 2$, 
then the assumption that $\mcM$ has the $\rho$-dimensional amalgamation property can be removed
and the conclusion still holds.
}\end{rem}

\section{Technical implications of trivial dependence in binary homogeneous structures}\label{implications of trivial dependence}

\noindent
In this section we define the notion of `$\acl$-complete set' and prove Theorem~\ref{homogeneity of D},
and its corollary, which shows, roughly speaking, that on any $\es$-definable $\acl$-complete subset of $\mcM\meq$ with rank 1 where $\mcM$ is
binary, homogeneous and simple with trivial dependence, the ``type structure'' is determined by the 2-types.

\begin{assump}\label{assumptions when working with D}{\rm
Throughout this section, including Theorem~\ref{homogeneity of D} and Corollary~\ref{set homogeneity of all of D},
we assume that
\begin{itemize}
\item[(i)] $\mcM$ is countable, binary, homogeneous, simple, with trivial dependence, and
\item[(ii)] $D \subseteq \mcM\meq$ is $\es$-definable, only finitely many sorts are represented in $D$,
and $\su(d) = 1$ for every $d \in D$.
\end{itemize}
}\end{assump}

\begin{notation}\label{notation in last two sections}{\rm
In the rest of the article, `$\tp_{\mcM\meq}$', `$\acl_{\mcM\meq}$' and `$\dcl_{\mcM\meq}$' are 
abbreviated with `$\tp$', `$\acl$' and `$\dcl$', respectively. (But when types, et cetera, are used with respect to other structures
we indicate it with a subscript.)
}\end{notation}

\noindent
Recall the notation `$F(\bar{a}, b)$' explained in the beginning of Section~\ref{Simple omega-categorical structures}.

\begin{lem}\label{types of imaginaries are determined by partial isomorphisms}
Suppose that $a_1, \ldots, a_n, b_1,  \ldots, b_n \in M\meq$.
Then the following are equivalent:
\begin{enumerate}
\item $\tp(a_1, \ldots, a_n) \ = \ \tp(b_1, \ldots, b_n)$.
\item There are finite sequences $\bar{a}_1, \ldots, \bar{a}_n, \bar{b}_1, \ldots, \bar{b}_n \in M$ and an isomorphism
$f : \mcM \uhrc \bar{a}_1 \ldots \bar{a}_n \to \mcM \uhrc \bar{b}_1 \ldots \bar{b}_n$ such that
$F(\bar{a}_i, a_i)$, $F(\bar{b}_i, b_i)$ and $f(\bar{a}_i) = \bar{b}_i$ for all $i = 1, \ldots, n$.
\end{enumerate}
\end{lem}

\noindent
{\em Proof.}
If $\tp(a_1, \ldots, a_n) \ = \ \tp(b_1, \ldots, b_n)$, then since $\mcM\meq$ is $\omega$-homogeneous and countable,
there is an automorphism $f$ of $\mcM\meq$ such that $f(a_i) = b_i$ for all $i$.
Let $\bar{a}_i \in M$ be such, for each $i$, that $F(\bar{a}_i, a_i)$, and let $\bar{b}_i = f(\bar{a}_i)$.
Then the restriction of $f$ to $\rng(\bar{a}_1) \cup \ldots \cup \rng(\bar{a}_n)$ is an isomorphism
from $\mcM \uhrc \bar{a}_1 \ldots \bar{a}_n$ to $\mcM \uhrc \bar{b}_1 \ldots \bar{b}_n$.

For the other direction, note that $F(\bar{a}_i, a_i)$ implies that $a_i \in \dcl(\bar{a}_i)$, and similarly for $\bar{b}_i$ and $b_i$.
So if~(2) holds then, as $\mcM$ is homogeneous, $\tp(\bar{a}_1, \ldots, \bar{a}_n) = \tp(\bar{b}_1, \ldots, \bar{b}_n)$, and therefore
$\tp(a_1, \ldots, a_n) = \tp(b_1, \ldots, b_n)$.
\hfill $\square$

\begin{defin}\label{definition of acl-completeness}{\rm
We call $D$ {\em $\acl$-complete} if for all $a \in D$ and all $\bar{a}, \bar{a}' \in M$, if 
$F(\bar{a}, a)$ and $F(\bar{a}', a)$, then $\tp(\bar{a} / \acl(a)) = \tp(\bar{a}' / \acl(a))$.
}\end{defin}

\begin{lem}\label{existence of acl-complete sets}
There is $D' \subseteq M\meq$ such that $D'$ satisfies Assumption~\ref{assumptions when working with D}~(ii), $D'$ is $\acl$-complete and
\begin{enumerate}
\item for every $d \in D$ there is (not necessarily unique) $d' \in D'$ such that $d \in \dcl(d')$ and $d' \in \acl(d)$, and
\item for every $d' \in D'$ there is $d \in D$ such that $d \in \dcl(d')$ and $d' \in \acl(d)$.
\end{enumerate}
\end{lem}

\noindent
{\em Proof.}
Let $p$ be any one of the finitely many complete 1-types over $\es$ which are realized in $D$,
and let the equivalence relation $E_p$ on $M^n$ (for some $n$) define the sort of the elements which realize $p$.
By Fact~\ref{second fact about canonically embedded structures}, the following equivalence relation on $M^n$ is $\es$-definable:
\[ E'_p(\bar{x}, \bar{y}) \ \Longleftrightarrow \ 
\exists z \Big( p(z) \ \wedge \ F_{E_p}(\bar{x}, z) \ \wedge \ F_{E_p}(\bar{y}, z)  \ \wedge \ \tp(\bar{x} / \acl(z)) = \tp(\bar{y} / \acl(z)) \Big).\]
Moreover, by the same fact, every $E_p$-class is a union of finitely many $E'_p$-classes.
By replacing, for every complete 1-type $p$ over $\es$ that is realized in $D$,
the elements realizing $p$ with the elements of $\mcM\meq$ which correspond to $E'_p$-classes,
we get $D'$. This set has the properties stated in the lemma.
\hfill $\square$
\\

\noindent
Recall Assumption~\ref{assumptions when working with D} and 
Notation~\ref{notation in last two sections}.
The following theorem is the main result of this section.

\begin{theor}\label{homogeneity of D}
Suppose that $D$ is $\acl$-complete, 
$1 < n < \omega$, $a_i, b_i \in D$ for $i = 1, \ldots, n$,
$\{a_1, \ldots, a_n\}$ is independent over $\es$,
$\{b_1,  \ldots, b_n\}$ is independent over $\es$ and
$\tp(a_i, a_j) = \tp(b_i, b_j)$ for all $i, j = 1, \ldots, n$.
Then  $\tp(a_1,  \ldots, a_n)  \ = \ \tp(b_1, \ldots, b_n)$.
\end{theor}

\begin{cor}\label{set homogeneity of all of D}
Suppose that $D$ is $\acl$-complete.
Let $n < \omega$,  $\bar{a}_i, \bar{b}_i \in D$, $i = 1, \ldots, n$, and suppose that
$\su(\bar{a}_i) = \su(\bar{b_i}) = 1$,
$\acl(\bar{a}_i) \cap D = \rng(\bar{a}_i)$
and $\acl(\bar{b}_i) \cap D = \rng(\bar{b}_i)$
for all $i$.
Furthermore, asssume that $\{\bar{a}_1, \ldots, \bar{a}_n\}$ is independent over $\es$,
$\{\bar{b}_1, \ldots, \bar{b}_n\}$ is independent over $\es$
and $\tp(\bar{a}_i, \bar{a}_j) = \tp(\bar{b}_i, \bar{b}_j)$ for all $i$ and $j$.
Then there is a permutation $\bar{b}'_i$ of $\bar{b}_i$, for each $i$, such that
$\tp(\bar{a}_1, \ldots, \bar{a}_n) = \tp(\bar{b}'_1, \ldots, \bar{b}'_n)$.
\end{cor}

\noindent
{\em Proof.}
Suppose that $\bar{a}_i = (a_{i,1}, \ldots, a_{i,k_i}), \bar{b}_i = (b_{i,1}, \ldots, b_{i,k_i})$, $i = 1, \ldots, n$, 
satisfy the assumptions of the theorem.
In particular we have $\tp(a_{i,1}, a_{j,1}) = \tp(b_{i,1}, b_{j,1})$ for all $i$ and $j$,
and both $\{a_{1,1}, \ldots, a_{n,1}\}$ and $\{b_{1,1}, \ldots, b_{n,1}\}$ are independent over $\es$.
By Theorem~\ref{homogeneity of D},
\[\tp(a_{1,1}, \ldots, a_{n,1}) \ = \ \tp(b_{1,1}, \ldots, b_{n,1}).\]
By $\omega$-homogeneity of $\mcM\meq$ (and Fact~\ref{equivalences to being homogeneous})
there is an automorphism $f$ of $\mcM\meq$ such that $f(a_{i,1}) = b_{i,1}$ for all $i = 1, \ldots, n$.
From $\su(\bar{a}_i) = 1$ it follows that $\bar{a}_i \in \acl(a_{i,1})$ for every $i$, and for the same reason
$\bar{b}_i \in \acl(b_{i,1})$ for every $i$.
Hence $f(\rng(\bar{a}_i)) = \rng(\bar{b}_i)$ for all $i$ and consequently
there is a permutation $\bar{b}'_i$ of $\bar{b}_i$ for each $i$ such that
$\tp(\bar{a}_1, \ldots, \bar{a}_n)  \ = \ \tp(\bar{b}'_1, \ldots, \bar{b}'_n)$. 
\hfill $\square$

\subsection{Proof of Theorem~\ref{homogeneity of D}}

Let $D \subseteq M\meq$ and $a_i, b_i \in D$, $i = 1, \ldots, n$, satisfy the assumptions of the theorem.
We prove that $\tp(a_1, \ldots, a_n) = \tp(b_1, \ldots, b_n)$ by induction on $n = 2, 3, 4, \ldots$. 
The case $n = 2$ are trivial, so we assume that $n > 2$ and, by the induction hypothesis, that
\begin{equation}\label{induction hypothesis for 2-determinedness of types}
\tp(a_1, \ldots a_{n-1}) = \tp(b_1, \ldots, b_{n-1}).
\end{equation}
Suppose that we can find $\bar{a}_i, \bar{b}_i \in M$, $i = 1, \ldots, n$, such that
$F(\bar{a}_i, a_i)$, $F(\bar{b}_i, b_i)$ for all $i$ and
\[\tp(\bar{a}_1, \ldots, \bar{a}_n) \ = \ 
\tp(\bar{b}_1, \ldots, \bar{b}_n).\]
Then
Lemma~\ref{types of imaginaries are determined by partial isomorphisms}
implies that 
\[\tp(a_1,  \ldots, a_n)  \ = \ \tp(b_1, \ldots, b_n)\]
which is what we want to prove.
Our aim is to find $\bar{a}_1,  \ldots, \bar{a}_n, \bar{b}_1, \ldots, \bar{b}_n$ as above.
We now prove three technical lemmas. Then a short argument which combines these lemmas proves the theorem.

\begin{lem}\label{first auxilliary lemma}
There are $\bar{a}_i \in M$ for $i = 1, \ldots, n$ such that $F(\bar{a}_i, a_i)$ for every $i$
and $\{\bar{a}_1, \ldots, \bar{a}_n\}$ is independent over $\es$.
\end{lem}

\noindent
{\em Proof.}
By induction we prove that for each $k = 1, \ldots, n$, there are $\bar{a}'_1, \ldots, \bar{a}'_k \in M$ such that
$F(\bar{a}'_i, a_i)$ for every $i$ and $\{\bar{a}'_1, \ldots, \bar{a}'_k\}$ is independent over $\es$.
The case $k = 1$ is trivial, so we assume that $0 < k < n$ and that we have found 
$\bar{a}_1, \ldots, \bar{a}_k \in M$ such that $F(\bar{a}_i, a_i)$ for every $i$ and
$\{\bar{a}_1, \ldots, \bar{a}_k\}$ is independent over $\es$.

Choose any $a^*_{k+1} \in D$ such that $\tp(a^*_{k+1} / a_1, \ldots, a_k, \bar{a}_1, \ldots, \bar{a}_k)$
is a nondividing extension of $\tp(a_{k+1} / a_1, \ldots, a_k)$, so in particular
$\tp(a^*_{k+1}, a_1, \ldots, a_k) = \tp(a_{k+1}, a_1, \ldots, a_k)$ and 
\[a^*_{k+1} \underset{a_1, \ldots, a_k}{\ind} \bar{a}_1, \ldots, \bar{a}_k\]
and, since (by assumption) $\{a_1, \ldots, a_{k+1}\}$ is independent over $\es$,
$a^*_{k+1} \ind a_1, \ldots, a_k$.
By transitivity of dividing,
\[a^*_{k+1} \ind a_1, \ldots, a_k, \bar{a}_1, \ldots, \bar{a}_k\] 
so by monotonicity
\begin{equation}\label{a-* is independent from the tuples}
a^*_{k+1} \ind \bar{a}_1, \ldots, \bar{a}_k.
\end{equation}
As $\tp(a^*_{k+1}, a_1, \ldots, a_k) = \tp(a_{k+1}, a_1, \ldots, a_k)$ and $\mcM\meq$ is $\omega$-homogeneous
(and countable) there is an automorphism $f$ of $\mcM\meq$ which maps $(a^*_{k+1}, a_1, \ldots, a_k)$ to $(a_{k+1}, a_1, \ldots, a_k)$.
Let $f(\bar{a}_i) = \bar{a}'_i$ for $i = 1, \ldots, k$.
Then $F(\bar{a}'_i, a_i)$ for $i = 1, \ldots, k$ and
\[\tp(a^*_{k+1}, a_1, \ldots, a_k, \bar{a}_1, \ldots, \bar{a}_k) \ = \ \tp(a_{k+1}, a_1, \ldots, a_k, \bar{a}'_1, \ldots, \bar{a}'_k),\]
so in view of~(\ref{a-* is independent from the tuples}),
\begin{equation}\label{a-k+1 is independent from the tuples}
a_{k+1} \ind \bar{a}'_1, \ldots, \bar{a}'_k,
\end{equation}
and as $\{\bar{a}_1, \ldots, \bar{a}_k\}$ is independent over $\es$ (by induction hypothesis),
\begin{equation}\label{the tuples with prime are independent}
\{\bar{a}'_1,  \ldots, \bar{a}'_k\} \ \text{ is independent over $\es$}.
\end{equation}

Choose any $\bar{a}_{k+1} \in M$ such that $F(\bar{a}_{k+1}, a_{k+1})$.
There are $\bar{a}^*_1,  \ldots, \bar{a}^*_k \in M\meq$ such that 
$\tp(\bar{a}^*_1, \ldots, \bar{a}^*_k / a_{k+1}, \bar{a}_{k+1})$ is a nondividing extension of
$\tp(\bar{a}'_1, \ldots, \bar{a}'_k / a_{k+1})$.
Then 
\[\bar{a}^*_1, \ldots, \bar{a}^*_k \underset{a_{k+1}}{\ind} \bar{a}_{k+1}\]
and $\tp(\bar{a}^*_1, \ldots, \bar{a}^*_k, a_{k+1}) = \tp(\bar{a}'_1, \ldots, \bar{a}'_k, a_{k+1})$,
so $\bar{a}^*_1, \ldots, \bar{a}^*_k \ind a_{k+1}$.
By transitivity, $\bar{a}^*_1, \ldots, \bar{a}^*_k  \ind a_{k+1}, \bar{a}_{k+1}$ and by monotonicity,
\begin{equation}\label{the bar-a-* are independent from a-k+1}
\bar{a}^*_1, \ldots, \bar{a}^*_k \ind \bar{a}_{k+1}.
\end{equation}
By the $\omega$-homogeneity of $\mcM\meq$ there is an automorphism $g$ of $\mcM\meq$ that maps
$(\bar{a}^*_1, \ldots, \bar{a}^*_k, a_{k+1})$ to $(\bar{a}'_1, \ldots, \bar{a}'_k, a_{k+1})$.
Let $g(\bar{a}_{k+1}) = \bar{a}'_{k+1}$.
Then 
\[\tp(\bar{a}^*_1, \ldots, \bar{a}^*_k,  \bar{a}_{k+1}, a_{k+1}) \ = \ 
\tp(\bar{a}'_1, \ldots, \bar{a}'_k, \bar{a}'_{k+1}, a_{k+1}),\]
so $F(\bar{a}'_{k+1}, a_{k+1})$ and, by~(\ref{the bar-a-* are independent from a-k+1}),
\[\bar{a}'_1, \ldots, \bar{a}'_{k} \ind \bar{a}'_{k+1}.\]
From~(\ref{the tuples with prime are independent}) it follows that
$\{\bar{a}'_1, \ldots, \bar{a}'_{k+1}\}$ is independent over $\es$.
\hfill $\square$

\begin{lem}\label{second auxilliary lemma}
Let $I \subseteq \{1, \ldots, n-1\}$.
Suppose that $\bar{c}_i  \in M$, for $i = 1, \ldots, n-1$, and $\bar{d}_j \in M$ for $j \in I$ 
are such that $F(\bar{c}_i, b_i)$ for every $1 \leq i \leq n-1$ and $F(\bar{d}_j, b_n)$ for every $j \in I$.
Then there are $\bar{c}'_i,  \in M$, for $i = 1, \ldots, n-1$, and $\bar{d}'_j \in M$ for $j \in I$ 
such that $F(\bar{c}'_i, b_i)$ for every $1 \leq i \leq n-1$, $F(\bar{d}'_j, b_n)$ for every $j \in I$,
$\tp(\bar{c}'_1, \ldots, \bar{c}'_{n-1}) = \tp(\bar{c}_1, \ldots, \bar{c}_{n-1})$,
$\tp(\bar{c}'_j, \bar{d}'_j) = \tp(\bar{c}_j, \bar{d}_j)$ for all $j \in I$ and
$b_n \notin \acl(\bar{c}'_1, \ldots, \bar{c}'_{n-1})$.
\end{lem}

\noindent
{\em Proof.}
Suppose on the contrary that $b_n \in \acl(\bar{c}'_1, \ldots, \bar{c}'_{n-1})$
for all $\bar{c}'_1, \ldots, \bar{c}'_{n-1} \in M$ such that 
\begin{align}\label{first equation in claim 2}
&\text{$F(\bar{c}'_i, b_i)$ for every $i = 1, \ldots, n-1$,} \\ 
&\tp(\bar{c}'_1, \ldots, \bar{c}'_{n-1}) = \tp(\bar{c}_1, \ldots, \bar{c}_{n-1}), \nonumber \\
&\text{for every $i \in I$ there is $\bar{d}'_i \in M$ such that $F(\bar{d}'_i, b_n)$, and} \nonumber \\
&\tp(\bar{c}'_i, \bar{d}'_i) = \tp(\bar{c}_i, \bar{d}_i). \nonumber
\end{align}
Note that by the $\omega$-categoricity of $\mcM$ the condition 
\[
\text{``$b_n \in \acl(\bar{c}'_1, \ldots, \bar{c}'_{n-1})$
for all $\bar{c}'_1, \ldots, \bar{c}'_{n-1} \in M$ such that~(\ref{first equation in claim 2}) holds''}
\] 
can be expressed by a
formula $\varphi(x_1,  \ldots, x_n)$ such that $\mcM\meq \models \varphi(b_1, \ldots, b_n)$.
By assumption, $\{b_1, \ldots, b_n\}$ is independent over $\es$, so $b_n \notin \acl(b_1, \ldots, b_{n-1})$ and hence there are
distinct $b_{n,i}$, for all $i < \omega$, such that
\[\tp(b_1,  \ldots, b_{n-1}, b_{n,i}) \ = \ \tp(b_1, \ldots, b_{n-1}, b_n) \ \text{ for all } i < \omega.\]
Then $\mcM\meq \models \varphi(b_1, \ldots, b_{n-1}, b_{n,i})$ for all $i < \omega$.
Since~(\ref{first equation in claim 2}) is satisfied if we let $\bar{c}'_i = \bar{c}_i$ for $i = 1, \ldots, n-1$ and
$\bar{d}'_i = \bar{d}_i$ for $i \in I$, it follows that
$b_{n,i} \in \acl(\bar{c}_1, \ldots, \bar{c}_{n-1})$ for all $i < \omega$.
This contradicts the $\omega$-categoricity of $\mcM$ 
(via Fact~\ref{only finitely many realizations}) 
because $\tp(b_{n,i}) = \tp(b_{n,j})$ for all $i$ and $j$.
\hfill $\square$
\\

\noindent
By Lemma~\ref{first auxilliary lemma}, let $\bar{a}_i \in M$ for $i = 1, \ldots, n$ be such that $F(\bar{a}_i, a_i)$ for every $i$
and 
\begin{equation}\label{the bar-a-primes are independent}
\{\bar{a}_1, \ldots, \bar{a}_n\} \ \text{ is independent over $\es$}.
\end{equation}

\begin{lem}\label{third auxilliary lemma}
Let $I$ be a proper subset of $\{1, \ldots, n-1\}$.
Suppose that $\bar{b}_i \in M$ for $i = 1, \ldots, n-1$ and $\bar{b}_{n,j} \in M$ for $j \in I$ are such that
\begin{align}\label{assumption of claim 3}
&F(\bar{b}_i, b_i) \ \text{ for all $i = 1, \ldots, n-1$}, F(\bar{b}_{n,j}, b_n) \ \text{ for all } j \in I, \\
&\tp(\bar{b}_1, \ldots, \bar{b}_{n-1}) \ = \ \tp(\bar{a}_1, \ldots, \bar{a}_{n-1}), \  \nonumber \\
&\tp(\bar{b}_j, \bar{b}_{n,j}) \ = \ \tp(\bar{a}_j, \bar{a}_n) \ \text{ for all } j \in I, \nonumber \\
&b_n \notin \acl(\bar{b}_1, \ldots, \bar{b}_{n-1}). \nonumber
\end{align}
Let $j \in \{1, \ldots, n-1\} \setminus I$ and $J = I \cup \{j\}$.
Then there are $\bar{b}'_i \in M$ for $i = 1, \ldots, n-1$ and $\bar{b}'_{n,j} \in M$ for $j \in J$ such that~(\ref{assumption of claim 3})
holds if `$\bar{b}$' is replaced with `$\bar{b}'$' and `$I$' with `$J$'.
\end{lem}

\noindent
{\em Proof.}
Suppose that~(\ref{assumption of claim 3}) holds. 
Note that the second line of it together with~(\ref{the bar-a-primes are independent}) 
implies that
\begin{equation}\label{the bar-b are independent over empty set}
\{\bar{b}_1, \ldots, \bar{b}_{n-1}\} \ \text{ is independent over $\es$}.
\end{equation}
Without loss of generality we assume that $I = \{1, \ldots, k\}$ where $k < n-1$. 
The case $k = 0$ is interpreted as meaning that $I = \es$.
By assumption
(of Theorem~\ref{homogeneity of D}), 
$\tp(b_{k+1}, b_n) = \tp(a_{k+1}, a_n)$, so there are $\bar{b}^*_{k+1}, \bar{b}_{n,k+1} \in M$ such that
\[\text{$F(\bar{b}^*_{k+1}, b_{k+1})$, $F(\bar{b}_{n,k+1}, b_n)$ and }
\tp(\bar{b}^*_{k+1}, \bar{b}_{n,k+1}, b_{k+1}, b_n) = \tp(\bar{a}_{k+1}, \bar{a}_n, a_{k+1}, a_n).\]
Since $\tp(\bar{b}_{n,k+1} / \bar{b}^*_{k+1}, b_n)$ has a nondividing extension to 
\[\bar{b}^*_{k+1}, b_n, \bar{b}_1, \ldots, \bar{b}_k, \bar{b}_{k+2}, \ldots, \bar{b}_{n-1}\]
we may without loss of generality assume that $\bar{b}_{n,k+1}$ realizes such a nondividing extension and hence
\begin{equation}\label{b-n,1 is independent over...}
\bar{b}_{n,k+1} \underset{b_n, \bar{b}^*_{k+1}}{\ind} \bar{b}_1, \ldots, \bar{b}_k, \bar{b}_{k+2}, \ldots, \bar{b}_{n-1}.
\end{equation}
From $\tp(\bar{b}^*_{k+1}, \bar{b}_{n,k+1}) = \tp(\bar{a}_{k+1}, \bar{a}_n)$ and~(\ref{the bar-a-primes are independent})
we get $\bar{b}^*_{k+1} \ind \bar{b}_{n,k+1}$, which since $b_{k+1} \in \dcl(\bar{b}^*_{k+1})$ implies that
$\bar{b}^*_{k+1}, b_{k+1} \ind \bar{b}_{n,k+1}$ and hence
\begin{equation}\label{independence of b-*-1 from b-n,1}
\bar{b}^*_{k+1} \underset{b_{k+1}}{\ind} \bar{b}_{n,k+1}.
\end{equation}
Since $b_{k+1} \in \dcl(\bar{b}_{k+1})$ it follows from~(\ref{the bar-b are independent over empty set}) that
\begin{equation}\label{independence of bar-b-1 from...}
\bar{b}_{k+1} \underset{b_{k+1}}{\ind} \bar{b}_1, \ldots, \bar{b}_k, \bar{b}_{k+2}, \ldots, \bar{b}_{n-1}.
\end{equation}
As $F(\bar{b}_{k+1}, b_{k+1})$ and $F(\bar{b}^*_{k+1}, b_{k+1})$, the assumption that $D$ is $\acl$-complete implies that
\begin{equation}\label{bar-b-1 and bar-b-*-1 have the same type over acl of b-1}
\tp(\bar{b}^*_{k+1} / \acl(b_{k+1})) = \tp(\bar{b}_{k+1} / \acl(b_{k+1})).
\end{equation}
We have already concluded that $\bar{b}_{n,k+1} \ind \, \bar{b}^*_{k+1}$ and since $b_n \in \dcl(\bar{b}_{n,k+1})$ we get
\[
\bar{b}_{n,k+1} \underset{b_n}{\ind} \bar{b}^*_{k+1},
\] 
which together with~(\ref{b-n,1 is independent over...}) and transitivity gives
\begin{equation}\label{bar-b-n,1 is independent from bar-b-*-1...}
\bar{b}_{n,k+1} \underset{b_n}{\ind} \bar{b}^*_{k+1}, \bar{b}_1, \ldots, \bar{b}_k, \bar{b}_{k+2}, \ldots, \bar{b}_{n-1}.
\end{equation}
Now we claim that
\begin{equation}\label{bar-b-n,1 is independent from bar-b-2}
\bar{b}_{n,k+1} \underset{b_{k+1}}{\ind} \bar{b}_1, \ldots, \bar{b}_k, \bar{b}_{k+2}, \ldots, \bar{b}_{n-1}.
\end{equation}
Suppose on the contrary that~(\ref{bar-b-n,1 is independent from bar-b-2}) is false.
Then $\bar{b}_{n,k+1} \nind \, b_{k+1}, \bar{b}_1, \ldots, \bar{b}_k, \bar{b}_{k+2}, \ldots, \bar{b}_{n-1}$.
Since $\bar{b}_{n,k+1} \ind \bar{b}^*_{k+1}$ (as we have seen above) and $b_{k+1} \in \dcl(\bar{b}^*_{k+1})$ we get $\bar{b}_{n,k+1} \ind b_{k+1}$.
By the triviality of dependence we must have $\bar{b}_{n,k+1} \nind \, \bar{b}_i$ for some $i \neq k+1$, 
so  
\[\bar{b}_{n,k+1}, b_n \nind \, \bar{b}_i.\]
Since $\su(b_n) = 1$ it follows from the last line of~(\ref{assumption of claim 3}) that $b_n \ind \bar{b}_i$.
From~(\ref{bar-b-n,1 is independent from bar-b-*-1...}) we get
$\bar{b}_{n,k+1} \underset{b_n}{\ind} \bar{b}_i$, so by transitivity
$\bar{b}_{n,k+1}, b_n \ind \bar{b}_i$ which contradicts what we got above.
Hence~(\ref{bar-b-n,1 is independent from bar-b-2}) is proved.

By the independence theorem
(Fact~\ref{independence theorem}) applied over $\acl(b_{k+1})$ together 
with (\ref{independence of b-*-1 from b-n,1}), (\ref{independence of bar-b-1 from...}),
(\ref{bar-b-1 and bar-b-*-1 have the same type over acl of b-1})
and~(\ref{bar-b-n,1 is independent from bar-b-2}), 
there is $\bar{b}'_{k+1}$ such that
\begin{align*}
&\tp(\bar{b}'_{k+1}, \bar{b}_{n,k+1}) \ = \ \tp(\bar{b}^*_{k+1}, \bar{b}_{n,k+1}) \ = \ \tp(\bar{a}_{k+1}, \bar{a}_n) \ \ \text{ and} \\
&\tp(\bar{b}_1, \ldots, \bar{b}_k, \bar{b}'_{k+1}, \bar{b}_{k+2}, \ldots, \bar{b}_{n-1}) \ = \
 \tp(\bar{b}_1, \ldots, \bar{b}_{n-1})
\ = \ \tp(\bar{a}_1, \ldots, \bar{a}_{n-1}).
\end{align*}
By applying Lemma~\ref{second auxilliary lemma} with $I = \{1, \ldots, k\}$,
$\bar{c}_i = \bar{b}_i$ for $i \in \{1, \ldots, n-1\} \setminus \{k+1\}$, $\bar{c}_{k+1} = \bar{b}'_{k+1}$
and $\bar{d}_i = \bar{b}_{n,i}$ for $i \in I$, we find
$\bar{b}'_i$ for $i \in \{1, \ldots, n-1\}$ and $\bar{b}'_{n,j}$ for $j \in J = I \cup \{k+1\}$ 
such that~(\ref{assumption of claim 3}) holds with `$\bar{b}'$',  and `$J$' in the place of 
`$\bar{b}$' and `$J$', respectively.
\hfill $\square$
\\

\noindent
Now we are ready to complete the proof of Theorem~\ref{homogeneity of D}.
By induction on $k = 1, \ldots, n-1$ and applying Lemma~\ref{third auxilliary lemma}
with $I = \{1, \ldots, k\}$ for $k < n-1$, we find $\bar{b}_1, \ldots, \bar{b}_{n-1} \in M$
and $\bar{b}_{n,1}, \ldots, \bar{b}_{n,n-1} \in M$ such that~(\ref{assumption of claim 3})
holds with $I = \{1, \ldots, n-1\}$.
With use of~(\ref{the bar-a-primes are independent}) it follows that
$\bar{b}_{n,i} \ind \bar{b}_i$ for all $i = 1, \ldots, n-1$ and since $b_n \in \dcl(\bar{b}_{n,i})$ we get
\begin{equation}\label{each bar-b-i is independent from bar-b-n-i over...}
\bar{b}_{n,i} \underset{b_n}{\ind} \bar{b}_i \ \ \text{ for all } i = 1, \ldots, n-1.
\end{equation}
Since $D$ is acl-complete we have
\begin{equation}\label{all bar-b-n-i have the same type of acl of b-n}
\tp(\bar{b}_{n,i} / \acl(b_n)) \ = \ \tp(\bar{b}_{n,j} / \acl(b_n)) \ \text{ for all } i,j = 1, \ldots, n-1.
\end{equation}
Moreover, we claim that
\begin{equation}\label{the bar-b-i are independent over b-n}
\{\bar{b}_1, \ldots, \bar{b}_{n-1}\} \ \ \text{ is independent over } \{b_n\}.
\end{equation}
Suppose on the contrary that~(\ref{the bar-b-i are independent over b-n}) is false.
By triviality of dependence, $\bar{b}_i \underset{b_n}{\nind} \bar{b}_j$ for some $i \neq j$,
and hence $\bar{b}_i \nind b_n \bar{b}_j$. By triviality of dependence again,
$\bar{b}_i \nind b_n$ or $\bar{b}_i \nind \bar{b}_j$.
But $\bar{b}_i \nind b_n$ implies $b_n \in \acl(\bar{b}_i)$ (since $\su(b_n) = 1$),
which contradicts the choice of $\bar{b}_1, \ldots, \bar{b}_{n-1}$.
And $\bar{b}_i \nind \bar{b}_j$ also contradicts the choice of $\bar{b}_1, \ldots, \bar{b}_{n-1}$
since $\tp(\bar{b}_i, \bar{b}_j) = \tp(\bar{a}_i, \bar{a}_j)$ where $\bar{a}_i \ind \bar{a}_j$.
Hence~(\ref{the bar-b-i are independent over b-n}) is proved.

The independence theorem (Corollary~\ref{corollary to independence theorem})
together with~(\ref{each bar-b-i is independent from bar-b-n-i over...}),
~(\ref{all bar-b-n-i have the same type of acl of b-n}) and~(\ref{the bar-b-i are independent over b-n}),
imply that there is $\bar{b}_n \in M$ such that $F(\bar{b}_n, b_n)$ and
$\tp(\bar{b}_n, \bar{b}_i) = \tp(\bar{b}_{n,i}, \bar{b}_i) = \tp(\bar{a}_n, \bar{a}_i)$ for all $i = 1, \ldots, n-1$.
Moreover, by the choice of $\bar{b}_1, \ldots, \bar{b}_{n-1}$, 
$\tp(\bar{b}_1, \ldots, \bar{b}_{n-1}) = \tp(\bar{a}_1, \ldots, \bar{a}_{n-1})$.
As the language is binary, there is an isomorphism $f$ from 
$\mcM \uhrc \bar{a}_1 \ldots \bar{a}_n$ to $\mcM \uhrc \bar{b}_1 \ldots \bar{b}_n$
such that $f(\bar{a}_i) = \bar{b}_i$ for each $i$,
so by Lemma~\ref{types of imaginaries are determined by partial isomorphisms},
$\tp(a_1, \ldots, a_n) = \tp(b_1, \ldots, b_n)$ and the proof of Theorem~\ref{homogeneity of D} is finished.

\section{Trivial dependence implies that any canonically embedded geometry is a reduct of a binary random structure}
\label{geometries}

\noindent
We use the conventions of Notation~\ref{notation in last two sections} throughout this section.

\begin{theor}\label{homogeneity of geometries}
Let $\mcM$ be countable, binary, homogeneous and simple with trivial dependence.
Suppose that $G \subseteq M\meq$ is $A$-definable where $A \subseteq M$ is finite, only finitely many sorts are represented in $G$, 
$\su(a/A) = 1$ and $\acl(\{a\} \cup A) \cap G = \{a\}$ for every $a \in G$.
Let $\mcG$ denote the canonically embedded structure in $\mcM\meq$ over $A$ with universe $G$.
Then $\mcG$ is a reduct of a binary random structure.
\end{theor}

\subsection{Proof of Theorem~\ref{homogeneity of geometries}}
Let $\mcM$, $G \subseteq M\meq$ and $A \subseteq M$ be as assumed in the theorem.
By Remark~\ref{remark about adding constants}, we may without loss of generality assume that $A = \es$,
implying that $G$ is $\es$-definable in $\mcM\meq$ and that $\mcG$ is a canonically embedded structure in $\mcM\meq$ over $\es$.
By Lemma~\ref{existence of acl-complete sets} applied to $G$, there is $D \subseteq M\meq$ with rank 1 such that
$D$ is $\es$-definable, acl-complete and 
\begin{align}\label{conditions on D}
&\text{for every $a \in G$ there is $d \in D$ such that $a \in \dcl(d)$ and $d \in \acl(a)$, and} \\
&\text{for every $d \in D$ there is $a \in G$ such that $a \in \dcl(d)$ and $d \in \acl(a)$.} \nonumber
\end{align}

\begin{rem}\label{extension property for D}
{\rm Observe that the independence theorem implies the following:}
Suppose that  $n < \omega$, $\{a_1, \ldots, a_n\} \subseteq D$ is independent over $\es$, 
$b_1, \ldots, b_n \in D$ and $b_i \ind a_i$ for all $i = 1, \ldots, n$ and 
$\tp(b_i / \acl(\es)) = \tp(b_j / \acl(\es))$ for all $i$ and $j$.
Then there is $b \in D$ such that $\tp(b / \acl(\es)) = \tp(b_i / \acl(\es))$ 
and $\tp(b, a_i) = \tp(b_i, a_i)$ for all $i = 1, \ldots, n$, and $b \ind \{a_1, \ldots, a_n\}$.
\end{rem}

\noindent
Let $p_1, \ldots, p_r$ be all complete 1-types over $\acl(\es)$ which are realized in $D$, and let $p_{r+1}, \ldots, p_s$ be all 
complete 2-types over $\es$
which are realized in $D$ and, for each $r < i \leq s$, have the property that if $p_i(a,b)$, then $a \neq b$ and $\{a, b\}$ is independent.
For each $i = 1, \ldots, s$, let $R_i$ be a relation symbol with arity $1$ if $i \leq r$ and otherwise with arity 2.
Let $V = \{R_1, \ldots, R_s\}$ and let $\mcD$ denote the $V$-structure with universe $D$ such that for every $\bar{a} \in D$,
$\mcD \models R_i(\bar{a})$ if and only if $\mcM\meq \models p_i(\bar{a})$.

Now define $\mbK$ to be the class of all {\em finite} $V$-structures $\mcN$ such that there is an embedding 
$f : \mcN \to \mcD$ such that $f(N)$ is an independent set.
Let $\mbP_2$ be the class of all $\mcN \in \mbK$ such that $|N| \leq 2$.
Recall the definition of $\mbR\mbP_2$ in Section~\ref{random structures}.

\begin{lem}\label{K is an amalgamation class}
$\mbK = \mbR\mbP_2$, where $\mbR\mbP_2$ has the hereditary property and the amalgamation property.
\end{lem}

\noindent
{\em Proof.}
$\mbP_2$ is clearly a 1-adequate class,  so (as observed in Section~\ref{random structures}) 
$\mbR\mbP_2$ has the hereditary property and amalgamation property.
We clearly have $\mbK \subseteq \mbR\mbP_2$, so it remains to prove that $\mbR\mbP_2 \subseteq \mbK$.
For this it suffices to show that if $\mcN \subset \mcN' \in \mbR\mbP_2$, $N' = N \cup \{a\}$
and $f : \mcN \to \mcD$ is an embedding such that $f(N)$ is independent, 
then there is an embedding $f' : \mcN' \to \mcD$ which extends $f$ and $f'(N')$ is independent.
But this follows immediately from Remark~\ref{extension property for D} together with the definitions of the involved structures.
\hfill $\square$
\\

\noindent
By Lemma~\ref{K is an amalgamation class}, $\mbK = \mbR\mbP_2$ has the hereditary property and the amalgamation property,
so let $\mcF$ be the Fra\"{i}ss\'{e} limit of $\mbK$. Hence $\mcF$ is homogeneous and a binary random structure.
Since $\mcF$ is the Fra\"{i}ss\'{e} limit of $\mbK$, it follows that if $\mcN \subseteq \mcN' \in \mbK$ and $f : \mcN \to \mcF$ is an embedding,
then there is an embedding $f' : \mcN' \to \mcF$ which extends $f$.
By using this together with the definition of $\mbK$ and Remark~\ref{extension property for D} it is straightforward
to prove, by a back and forth argument, that there is $D' \subseteq D$ such that 
\begin{itemize}
\item[(a)] $D'$ is independent,
\item[(b)] $\mcF \cong \mcD \uhrc D'$, and
\item[(c)] for every $d \in D$ there is $d' \in D'$ such that $\acl(d) = \acl(d')$.
\end{itemize}
Let $a \in G$.
By~(\ref{conditions on D}), $a \in \dcl(d)$ for some $d \in D$.
By~(c), there is $d' \in D'$ such that $\acl(d) = \acl(d')$ and hence $a \in \acl(d')$.
For a contradiction suppose that there is $a' \in G$ such that $a' \neq a$, $\tp(a') = \tp(a)$ and $a' \in \acl(d')$.
By~(\ref{conditions on D}) and~(c) there is $d'' \in D' \cap \acl(a')$.
As $a' \in \acl(d')$ this implies that  $d'' \in \acl(d')$, which by the independence of $D'$ gives $d'' = d'$.
Then $a \in \acl(d') =  \acl(d'') \subseteq  \acl(a')$ which contradicts the assumptions about $G$.
Thus we conclude that
\begin{itemize}
\item[(d)] every $a \in G$ belongs to $\dcl(d')$ for some $d' \in D'$.
\end{itemize}
From the assumptions about $G$, $D$ and~(a) it follows that for every $a \in G$, $\acl(a) \cap D'$ contains a unique element
which we denote $g(a)$. It also follows from the assumptions about $G$, $D$ and~(a) that 
$g : G \to D'$ is bijective, and, using~(d), that
\begin{itemize}
\item[(e)] for every $a \in G$, $a \in \dcl(g(a))$.
\end{itemize}
Observe that we are not assuming, and we have not proved, that $D'$ or $g$ are definable (over any set).

Now we define a $V$-structure $\mcG'$ with universe $G$ as follows. 
For each $R_i \in V$ and every $\bar{a} \in G$, let 
\begin{itemize}
\item[] $\mcG' \models R_i(\bar{a})$ if and only if $\mcD \uhrc D' \models R_i(g(\bar{a}))$.
\end{itemize} 
Since $g$ is bijective it is clear that $\mcG' \cong \mcD \uhrc D'$ and by~(b) we get $\mcG' \cong \mcF$ so 
$\mcG'$ is a binary random structure.
From the definition of $\mcG'$ (through the definitions of $\mcD$, $D'$ and $\mcF$) it follows that
for every $a \in G$ there is $R_i$, $1 \leq i \leq r$, such that $\mcG \models R_i(a)$, and for all distinct $a, b \in G$
there is $R_i$, $r < i \leq s$, such that $\mcG \models R_i(a,b)$.

\begin{lem}\label{types in G are determined by types in G'}
If $n < \omega$, $a_1, \ldots, a_n, b_1, \ldots, b_n \in G$ and 
$\tp_{\mcG'}(a_1, \ldots, a_n) = \tp_{\mcG'}(b_1, \ldots, b_n)$, then $\tp_\mcG(a_1, \ldots, a_n) = \tp_\mcG(b_1, \ldots, b_n)$.
\end{lem}

\noindent
{\em Proof.}
Suppose that $a_1, \ldots, a_n, b_1, \ldots, b_n \in G$ and $\tp_{\mcG'}(a_1, \ldots, a_n) = \tp_{\mcG'}(b_1, \ldots, b_n)$.
Since $\mcG$ is a canonically embedded structure in  $\mcM\meq$, it follows that 
$\tp_\mcG(\bar{a}) = \tp_\mcG(\bar{b})$ if and only if $\tp(\bar{a}) = \tp(\bar{b})$, for all finite tuples $\bar{a}, \bar{b} \in G$.
So it suffices to prove that $\tp(a_1, \ldots, a_n) = \tp(b_1, \ldots, b_n)$.
We may assume that all $a_1, \ldots, a_n$ are distinct and that all $b_1, \ldots, b_n$ are distinct.

The assumptions and the definitions of $\mcG$, $\mcD$ and $D'$ imply that 
$\tp(g(a_i), g(a_j)) = \tp(g(b_i), g(b_j))$ for all $i$ and $j$.
Since $g : G \to D'$ is bijective and $D'$ is independent it follows from 
Theorem~\ref{homogeneity of D} that
\[\tp(g(a_1), \ldots, g(a_n)) \ = \ \tp(g(b_1), \ldots, g(b_n)).\]
By~(e) we have $a_i \in \dcl(g(a_i))$ and $b_i \in \dcl(g(b_i))$ for each $i$, and therefore 
$\tp(a_1, \ldots, a_n) = \tp(b_1, \ldots, b_n)$.
\hfill $\square$
\\

\noindent
To prove that $\mcG$ is a reduct of $\mcG'$ it suffices to show that for every $1 < n < \omega$ and every complete
$n$-type over $\es$ of $\mcG$ there is a $V$-formula $\varphi_p(\bar{x})$ such that for all $n$-tuples $\bar{a} \in G$,
$\mcG \models p(\bar{a})$ if and only if $\mcG' \models \varphi_p(\bar{a})$.
As $\mcG'$ has elimination of quantifiers it has only finitely many complete $n$-types over $\es$, say $q_1, \ldots, q_m$.
Let $q_i$ be isolated by $\varphi_i(\bar{x})$.
By Lemma~\ref{types in G are determined by types in G'},
for each $i$ either
\begin{itemize}
\item for all $\bar{a} \in G$, if $\mcG' \models \varphi_i(\bar{a})$, then $\mcG \models p(\bar{a})$, or
\item for all $\bar{a} \in G$, if $\mcG' \models \varphi_i(\bar{a})$, then $\mcG \not\models p(\bar{a})$
\end{itemize}
Let $I$ be the set of all $i$ for which the first case holds.
If $\varphi_p(\bar{x}) = \bigvee_{i \in I} \varphi_i(\bar{x})$ then, for all $n$-tuples $\bar{a} \in G$,
$\mcG \models p(\bar{a})$ if and only if $\mcG' \models \varphi_p(\bar{a})$.
This concludes the proof of Theorem~\ref{homogeneity of geometries}.

\begin{rem}\label{when is D' definable}{\rm
The conclusion of Theorem~\ref{homogeneity of geometries} is that 
\begin{itemize}
\item[(f)] $\mcG$ is a reduct of a binary random structure.
\end{itemize}
A stronger conclusion, essentially saying that $\mcG$ is a binary random structure would be:
\begin{itemize}
\item[(g)] If $\mcG_0$ is the reduct of $\mcG$ to all relation symbols with arity at most 2, then $\mcG$ is a reduct of $\mcG_0$,
and $\mcG_0$ is a binary random structure.
\end{itemize}
What extra assumptions do we need in order to get the conclusion~(g)?
It is straightforward to verify the following implications, where we use notation from the above proof:
\begin{align*}
&\text{ $D'$ is $\es$-definable in $\mcM\meq$} \\
\Longleftrightarrow 
&\text{ (the graph of) $g$ is $\es$-definable} \\
\Longrightarrow
&\text{ $g(a) \in \dcl(a)$ for every $a \in G$} \\
\Longrightarrow
&\text{ for all $0 < n < \omega$ and all $a_1, \ldots, a_n, b_1, \ldots, b_n \in G$,} \\
&\text{ $\tp(a_1, \ldots, a_n) = \tp(b_1, \ldots, b_n)$ if and only if } \\
&\text{ $\tp(g(a_1), \ldots, g(a_n)) = \tp(g(b_1), \ldots, g(b_n))$.}
\end{align*}
It follows that the condition that $D'$ is definable over $\es$ in $\mcM\meq$, as well as the equivalent condition, guarantees
that the conclusion of the proof of Theorem~\ref{homogeneity of geometries} is~(g).
The next example shows that~(g) does not in general follow from the assumptions of
Theorem~\ref{homogeneity of geometries}.
}\end{rem}

\begin{exam}\label{example of the referee}{\rm
This example, due to  the anonymous referee, shows that there are $\mcM$ and $G \subseteq M\meq$
which satisfy the assumptions of Theorem~\ref{homogeneity of geometries} but for which (g) fails if 
we let $\mcG$ be the canonically embedded structure (in $\mcM\meq$) with universe $G$.
Consequently, for such  $\mcM$ and $G$ every $D'$ as in the proof
of Theorem~\ref{homogeneity of geometries} is {\em not} $\es$-definable.

Let $\mcF$ be the random graph, that is, $\mcF$ is the Fra\"{i}ss\'{e} limit of the class of all finite
undirected loopless graphs. We now construct a new graph $\mcM$ (viewed as a first-order structure) as follows.
The universe of $\mcM$ is $M = F \times \{0,1\}$ (where $F$ is the universe of $\mcF$).
If $a, b \in F$ are adjacent then $(a, i)$ and $(b, i)$ are adjacent in $\mcM$ for $i = 0, 1$.
If $a, b \in F$ are nonadjacent (so in particular if $a = b$) then $(a, i)$ and $(b, 1-i)$ are adjacent in $\mcM$ for $i = 0, 1$.
There are no other adjacencies in $\mcM$.

Now we define
\begin{itemize}
\item[] for $(a,i), (b,j) \in M$, $E((a,i), (b,j))$ if and only if $a = b$.
\end{itemize}
Clearly $E$ is an equivalence relation such that each one of its classes has cardinality 2.
Moreover, it is straightforward to see that $E(x,y)$ is $\es$-definable by the formula
\[
x = y \ \vee \ \big( x \neq y \ \wedge \ \neg \exists z \big( z \sim_\mcM x \ \wedge \ z \sim_\mcM y \big) \big),
\]
where `$\sim_\mcM$' denotes adjacency in $\mcM$.
For every $u \in M$ let $u'$ denote the unique $v \neq u$ such that $E(u, v)$ holds
(or in other words, for $u=(a, i) \in M$, $u' = (a, 1-i)$).
Note that for all $u \in M$, $(u')' = u$ and $u \sim_\mcM u'$.
Let 
\[
M_0 \ = \ \{ (a,0) : a \in F\} \ \ \text{ and } \ \ M_1 \ = \ \{ (a,1) : a \in F\}
\]
and note that the set $M$ is the disjoint union of $M_0$ and $M_1$ and that $\mcM \uhrc M_i$ is a copy
of the random graph for $i = 0,1$.
The following is a straightforward consequence of the definition of  $\mcM$:

\medskip
\noindent
{\em Claim 1: For all distinct $u, v \in M$, 
\[
u \sim_\mcM v \ \Longleftrightarrow \ u' \sim_\mcM v' \ \Longleftrightarrow \
u \not\sim_\mcM v' \ \Longleftrightarrow \ u' \not\sim_\mcM v.
\]
}

\medskip

We now prove that $\mcM$ is homogeneous.
The above claim tells that  if $n < \omega$, $u_1, \ldots, u_n$, $v_1, \ldots, v_n \in M$ and
$f(u_i) = v_i$ for $i = 1, \ldots, n$ is a partial isomorphism, then $f$ can be extended
to a partial isomorphism which maps $u'_i$ to $v'_i$ for all $i = 1, \ldots, n$.
So to prove that $\mcM$ is homogeneous it suffices (by the symmetry of $M_0$ and $M_1$) to prove the following:

\medskip

\noindent
{\em Claim 2: Let $u_1, \ldots u_n, v_1, \ldots, v_n \in M_0$ and suppose that the map
$f(u_i) = v_i$ and $f(u'_i) = v'_i$ for $i = 1, \ldots, n$ is a partial isomorphism. Then for every
$u_{n+1} \in M_0$ there is $v_{n+1} \in M_0$ such that $f$ can be extended to a partial isomorphism
$g$ such that $g(u_{n+1}) = v_{n+1}$ and $g(u'_{n+1}) = g(v'_{n+1})$.}

\medskip
\noindent
We do not give the details of the proof of this claim but just note that the argument is straightforward and
uses that $\mcF$ is the random graph, the construction of $\mcM$ and the first claim.

By representing 0 and 1 with two distinct elements of $a_0, a_1 \in F$ it is straightforward to verify that
$\mcM$ is interpretable in $\mcF$ with the parameters $a_0$ and $a_1$.
It follows (from \cite[Remarks 2.26 and 2.27]{Cas}) that $\mcM$ is simple.
The natural way of interpreting $\mcM$ in $\mcF$ (with the parameters $a_0, a_1$) is by
letting $F^- = F \setminus \{a_0, a_1\}$, so $\mcF \uhrc F^- \cong \mcF$, and then identifying
the universe of $M$ with $F^- \times \{a_0, a_1\}$. Then $\su(u / \{a_0, a_1\}) = 1$ for every $u \in M$
(where SU-rank is with respect to $\mcF$)
and it follows that $\mcM$ is supersimple with SU-rank 1. 
Moreover, since $\mcF$ has trivial dependence it follows that the same is true for $\mcM$.
Because if there where subsets of $N\meq$ for some $\mcN \equiv \mcM$ that witnessed nontrivial dependence, then, 
by supersimplicity, we may assume that they are finite, so by $\omega$-categoricity of $\mcM$ we may assume
that they are subsets of $\mcM\meq$, and finally
the same sets with $a_0$ and $a_1$ added would witness nontrivial dependence in $\mcF$, a contradiction.
Hence $\mcM$ satisfies the assumptions of Theorem~\ref{homogeneity of geometries}.

For every $u \in M$ let $[u]$ be its equivalence class with respect to $E$.
Let 
\[ G = \{[u] : u \in M\}.\]
Then $G \subseteq M\meq$ and $G$ satisfies the assumptions of Theorem~\ref{homogeneity of geometries}.
Let $\mcG$ be the canonically embedded structure with universe $G$ and let 
$\mcG_0$ be the reduct of $\mcG$ to the relation symbols of arity at most 2.
It remains prove that $\mcG$ is {\em not} a reduct of $\mcG_0$.

First we show the following:

\medskip
\noindent
{\em Claim 3: For all distinct $u_1, u_2 \in G$ and all distinct $v_1, v_2 \in G$,
$\tp_\mcG(u_1, u_2) = \tp_\mcG(v_1, v_2)$.}

\medskip
\noindent
Let  $g_1, g_2, h_1, h_2 \in G$ be such that $g_1 \neq g_2$ and $h_1 \neq h_2$.
Then there are $u_1, u_2, v_1, v_2 \in M_0$ such that
$g_i = \{ u_i, u'_i \}$ and $h_i = \{v_i, v'_i\}$ for $i = 1,2$.

We consider four cases: (1) $u_1 \sim_\mcM u_2$ and $v_1 \sim_\mcM v_2$, 
(2) $u_1 \not\sim_\mcM u_2$ and $v_1 \not\sim_\mcM v_2$,
(3) $u_1 \sim_\mcM u_2$ and $v_1 \not\sim_\mcM v_2$, and
(4) $u_1 \not\sim_\mcM u_2$ and $v_1 \sim_\mcM v_2$.
In the first two cases the map given by $u_i \mapsto v_i$ and $u'_i \mapsto v'_i$ for $i = 1,2$ is a partial 
isomorphism, so by the homogeneity of $\mcM$ it extends to an automorphism of $\mcM$ and hence we get 
$\tp_\mcM(u_1, u_2, u'_1, u'_2) = \tp_\mcM(v_1, v_2, v'_1, v'_2)$ which in turn gives
$\tp_\mcG(g_1, g_2) = \tp_\mcG(h_1, h_2)$ (since $g_i \in \dcl_{\mcM\meq}(u_i)$ and similarly for $h_i$).
In the third and fourth case the map given by
$u_1 \mapsto v_1$, $u'_1 \mapsto v'_1$, $u_2 \mapsto v'_2$ and $u'_2 \mapsto v_2$ is a partial isomorphism
so we get 
$\tp_\mcM(u_1, u_2, u'_1, u'_2) = \tp_\mcM(v_1, v'_2, v'_1, v_2)$ and hence
$\tp_\mcG(g_1, g_2) = \tp_\mcG(h_1, h_2)$ (as $h_2 \in \dcl_{\mcM\meq}(v'_2)$).
This concludes the proof of Claim~3.

Observe that Claim~3 implies that every isomorphism between finite substructures
of $\mcG_0$ can be extended to an automorphism of $\mcG_0$, 
so $\mcG_0$ is a binary homogeneous structure (with finite vocabulary). 
This and Claim~3 easily implies the following:

\medskip
\noindent
{\em Claim 4: For every $n < \omega$, all distinct $g_1, \ldots, g_n \in G$ and
all distinct $h_1, \ldots, h_n \in G$,
$\tp_{\mcG_0}(g_1, \ldots, g_n) = \tp_{\mcG_0}(h_1, \ldots, h_n)$.}

\medskip
To prove that $\mcG$ is not a reduct of $\mcG_0$ 
it now suffices to show that there are distinct $g_1, g_2, g_3 \in G$ and distinct $h_1, h_2, h_3 \in G$
such that $\tp_\mcG(g_1, g_2, g_3) \neq \tp_\mcG(h_1, h_2, h_3)$.
Since $\mcM$ restricted to $M_0$ is a copy of the random graph it follows that there are distinct
$u_1, u_2, u_3 \in M_0$ and distinct $v_1, v_2,  v_3 \in M_0$ 
such that 
\[
u_1 \sim_\mcM u_2, \ u_1 \sim_\mcM u_3, \ u_2 \not\sim_\mcM u_3 \ \ 
\text{and $v_1, v_2, v_3$ forms a 3-cycle}.
\]
By Claim~1 we see that 
\[
\mcM \uhrc \{u_1, u_2, u_3, u'_1, u'_2, u'_3\} \not\cong \mcM \uhrc \{v_1, v_2, v_3, v'_1, v'_2, v'_3\}.
\]
Let $g_i = [u_i]$ and $h_i = [v_i]$ for $i = 1,2,3$.
Then $g_1, g_2, g_3$ are distinct and the same holds for $h_1, h_2, h_3$. 
Moreover,
\begin{align*}
&\acl_{\mcM\meq}(g_1, g_2, g_3) \cap M \ = \ \{u_1, u_2, u_3, u'_1, u'_2, u'_3\} \ \ \text{ and} \\
&\acl_{\mcM\meq}(h_1, h_2, h_3) \cap M \ = \ \{v_1, v_2, v_3, v'_1, v'_2, v'_3\}.
\end{align*}
It follows that $\tp_\mcG(g_1, g_2, g_3) \neq \tp_\mcG(h_1, h_2, h_3)$
and this finishes the proof that this example has the claimed properties.

We know from Theorem~\ref{homogeneity of geometries} that $\mcG$ is
a reduct of a binary random structure. In this example we can explicitly describe such a binary random structure.
We can simply expand $\mcM$ with a unary relation symbol interpreted as $M_0$.
Call this expansion $\mcM^*$. 
Let $\mcG^*$ be the canonically embedded structure of $(\mcM^*)\meq$ with universe $G$.
Let $\mcG^*_0$ be the reduct of $\mcG^*$ to the relation symbols of arity at most $2$.
One can now prove that $\mcG^*_0$ is a binary random structure and that 
$\mcG$ is a reduct of $\mcG^*_0$. 
}\end{exam}

\noindent
{\bf Acknowledgement.} 
We thank the anonymous referee for supplying Example~\ref{example of the referee} 
and for careful reading of the article.


\begin{thebibliography}{99}\label{References}

\bibitem{AL} A. Aranda L\'{o}pez, {\em Omega-categorical simple theories}, Ph.D. thesis, 
The University of Leeds (2014).

\bibitem{Cas} E. Casanovas, {\em Simple theories and hyperimaginaries}, Lecture Notes in Logic 39, 
The Association for Symbolic Logic and Cambridge University Press (2011)

\bibitem{Che98} G. L. Cherlin, {\em The Classification of Countable Homogeneous Directed Graphs
and Countable Homogeneous $n$-tournaments}, 
Memoirs of the American Mathematical Society 621, American Mathematical Society (1998).

\bibitem{CH} G. Cherlin, E. Hrushovski, {\em Finite Structures with Few Types}, 
Annals of Mathematics Studies 152, Princeton University Press (2003).

\bibitem{CHL} G. Cherlin, L. Harrington, A. H. Lachlan, {\em $\aleph_0$-categorical, $\aleph_0$-stable structures},
Annals of Pure and Applied Logic, Vol. 28 (1985) 103--135.

\bibitem{PK} T. De Piro, B. Kim, {\em The geometry of 1-based minimal types},
Transactions of The American Mathematical Society, Vol. 355 (2003) 4241--4263.

\bibitem{Djo06} M. Djordjevi\'{c}, {\em Finite satisfiability and $\omega$-categorical structures with trivial dependence},
 The Journal of Symbolic Logic, Vol. 71 (2006) 810--829.

\bibitem{EF} H-D. Ebbinghaus, J. Flum, {\em Finite Model Theory}, Second Edition, Springer-Verlag (1999).

\bibitem{Fra54} R. Fra\"{i}ss\'{e}, {\em Sur l'extension aux relations de quelques propri\'{e}t\'{e}s des ordres},
Annales Scientifiques de l'\'{E}cole Normale Sup\'{e}rieure, Vol. 71 (1954) 363--388.

\bibitem{Gar} A. Gardiner, {\em Homogeneous graphs}, 
Journal of Combinatorial Theory, Series B, Vol. 20 (1976) 94--102.

\bibitem{GK} Y. Golfand, M. Klin, {\em On $k$-homogeneous graphs}, 
in Algorithmic Studies in Combinatorics, Nauka, Moscow (1978), 76--85.

\bibitem{HKP} B. Hart, B. Kim, A. Pillay, {\em Coordinatisation and canonical bases in simple theories},
The Journal of Symbolic Logic, Vol. 65 (2000) 293--309.

\bibitem{Hen72} C. W. Henson, {\em Countable homogeneous relational structures and $\omega$-categorical theories},
The Journal of Symbolic Logic, Vol. 37 (1972) 494--500.

\bibitem{Hod} W. Hodges,  {\em Model theory}, Cambridge University Press (1993).

\bibitem{JTS} T. Jenkinson, J. K. Truss, D. Seidel, {\em Countable homogeneous multipartite graphs}, 
European Journal of Combinatorics, Vol. 33 (2012) 82--109. 

\bibitem{KLM} W. M. Kantor, M. W. Liebeck, H. D. Macpherson, 
{\em $\aleph_0$-categorical structures smoothly approximated by finite structures},
Proceedings of the London Mathematical Society, Vol. 59 (1989) 439--463.

\bibitem{Kol05} A. S. Kolesnikov, {\em $n$-Simple theories},  Annals of Pure and Applied Logic,
Vol. 131 (2005) 227--261.

\bibitem{Kop12} V. Koponen, {\em Asymptotic probabilities of extension properties and random $l$-colourable structures},
Annals of Pure and Applied Logic, Vol. 163 (2012) 391--438.

\bibitem{Kop-one-based} V. Koponen, {\em Homogeneous 1-based structures and interpretability in random structures}, submitted.

\bibitem{Lach84} A. H. Lachlan,  {\em Countable homogeneous tournaments},
Transactions of the American Mathematical Society, Vol. 284 (1984) 431--461.

\bibitem{Lach97} A. H. Lachlan, {\em Stable finitely homogeneous structures: a survey},
in B. T. Hart et. al. (eds.), Algebraic Model Theory, 145--159, Kluwer Academic Publishers (1997) 

\bibitem{LT} A. H. Lachlan, A. Tripp, {\em Finite homogeneous 3-graphs}, 
Mathematical Logic Quarterly, Vol. 41 (1995) 287--306.

\bibitem{LW} A. H. Lachlan, R. Woodrow, {\em Countable ultrahomogenous undirected graphs},
Transactions of the Americal Mathematical Society, Vol. 262 (1980) 51--94.

\bibitem{Mac91} D. Macpherson, {\em Interpreting groups in $\omega$-categorical structures},
The Journal of Symbolic Logic, Vol. 56 (1991) 1317--1324.

\bibitem{Mac10} D. Macpherson, {\em  A survey of homogeneous structures},
Discrete Mathematics, Vol. 311 (2011) 1599--1634.

\bibitem{Ober} W. Oberschelp, {\em Asymptotic 0-1 laws in combinatorics}, in D. Jungnickel (Ed.),
Combinatorial Theory, Lecture Notes in Mathematics, Vol. 969, Springer-Verlag (1982) 276--292.

\bibitem{Schm} J. H. Schmerl, {\em Countable homogeneous partially ordered sets},
Algebra Universalis, Vol. 9 (1979) 317--321.

\bibitem{Shee} J. Sheehan, {\em Smoothly embeddable subgraphs},
Journal of The London Mathematical Society, Vol. 9 (1974) 212--218.

\bibitem{She} S. Shelah, {\em Classification Theory}, Revised Edition, North-Holland (1990).

\bibitem{Wag} F. O. Wagner, {\em Simple Theories}, Kluwer Academic Publishers (2000).

\end{thebibliography}
\end{document}